\documentclass{article}    
\usepackage{amssymb, latexsym, amsmath, amsthm}

\newcommand{\bi} {\begin{itemize}}
\newcommand{\ei} {\end{itemize}}
\newcommand{\bea} {\begin{eqnarray}}
\newcommand{\eea} {\end{eqnarray}}
\newcommand{\be} {\begin{equation}}
\newcommand{\ee} {\end{equation}}
\newcommand{\bean} {\begin{eqnarray*}}
\newcommand{\eean} {\end{eqnarray*}}

\RequirePackage{epsfig}
\RequirePackage{amssymb}
\RequirePackage{natbib}
\RequirePackage{graphics}
\bibpunct{(}{)}{,}{a}{,}{,}

\begin{document}
\title{Boundary values in $R^t(K,\mu)$-spaces and invariant subspaces}
\author{Liming Yang}
\author{
Liming Yang \\
Department of Mathematics \\
Virginia Polytechnic Institute and State University \\
Blacksburg, VA 24061 \\
yliming@vt.edu
     }         
\date{}
\maketitle

\newtheorem{Theorem}{Theorem}
\newtheorem*{MTheorem}{Main Theorem}
\newtheorem*{TTheorem}{Thomson's Theorem}
\newtheorem*{ARSTheorem}{Aleman-Richter-Sundberg's Theorem}
\newtheorem*{Conjecture}{Conjecture}
\newtheorem{Corollary}{Corollary}
\newtheorem{Definition}{Definition}
\newtheorem{Assumption}{Assumption}
\newtheorem{Lemma}{Lemma}
\newtheorem{Problem}{Problem}
\newtheorem*{Example}{Example}
\newtheorem*{Remark}{Remark}
\newtheorem{KnownResult}{Known Result}
\newtheorem{Algorithm}{Algorithm}
\newtheorem{Property}{Property}
\newtheorem{Proposition}{Proposition}


\abstract{
For $1 \le t < \infty ,$ a compact subset $K$ of the complex plane $\mathbb C,$ and a finite positive measure $\mu$ supported on $K,$ $R^t(K, \mu)$ denotes the closure in $L^t (\mu )$ of rational functions with poles off $K.$ The paper examines the boundary values of functions in $R^t(K, \mu)$ for certain compact subset $K$ and extends the work of Aleman, Richter, and Sundberg on nontangential limits for the closure in $L^t (\mu )$ of analytic polynomials (Theorem A and Theorem C in \cite{ars}). We show that the Cauchy transform of an
annihilating measure has some continuity properties in the sense of capacitary density. This allows us to extend Aleman, Richter, and Sundberg's results for $R^t(K, \mu)$ and provide alternative short proofs of their theorems as special cases.
 
}

\section{Introduction}

Let $\mathcal{P}$ denote the set of polynomials in the complex variable $z.$ For a compact subset $K$ of the complex plane $\mathbb C,$ let $Rat(K)$ be the set of rational functions with poles off $K.$
For $1 \le t < \infty$ with conjugate exponent $t' = \dfrac{t}{t-1}$ and a finite positive measure $\mu$ supported on $K,$ let $R^t(K, \mu)$ denote the closure in $L^t (\mu )$ of $Rat(K).$ In the case that $K$ is polynomially convex, $R^t(K, \mu) = P^t(\mu ), $ the closure of $\mathcal{P}$ in $L^t(\mu ).$ Multiplication by $z$ defines a bounded linear operator on $R^t(K, \mu)$ which we will denote by $S_\mu.$ A rationally  invariant subspace of $R^t(K, \mu)$ is a closed linear subspace $M \subset R^t(K, \mu)$ such that $r(S_\mu) M \subset M$ for $r\in Rat(K).$
For a subset $A \subset \mathbb C,$ we set $\bar A$ or $clos (A)$  for its closure, $A^c$ for its complement, and $\chi _A$ for its characteristic function. For $\lambda\in \mathbb C$ and $\delta > 0,$ we set $B(\lambda, \delta) = \{z: |z - \lambda  | <\delta \}$ and $\mathbb D = B(0,1).$ Let $m$ be the normalized Lebesgue measure $\frac{d\theta}{2\pi}$ on $\partial {\mathbb D}.$ For a compactly supported finite measure $\nu$ on $\mathbb {C},$ we denote the support of $\nu$ by $spt(\nu ).$ For a compact subset $K,$ we denote the boundary of $K$ by $\partial K.$ The inner boundary of $K$, denoted by $\partial _i K$, is the set of boundary points
which do not belong to the boundary of any connected component of $\mathbb C\setminus K.$

For $\lambda \in K,$ we denote evaluation on $Rat(K)$ at $\lambda$ by $e_\lambda ,$ i.e. $e_\lambda (r) = r(\lambda )$ for $r\in Rat(K).$ 
$\lambda$ is a bounded point evaluation (bpe) for $R^t(K, \mu)$ if $e_\lambda$  extends to a bounded linear functional on $R^t(K, \mu),$ which we will also denote by $e_\lambda .$ We denote the set of bounded point evaluations for $R^t(K, \mu)$ by $bpe(R^t(K, \mu))$ and set $M_\lambda = \|e_\lambda\|_{R^t(K, \mu)^*} .$
For $\lambda _0 \in K,$ if there is a neighborhood of $\lambda _0,$ $B(\lambda _0, \delta),$ consisting entirely
of bpe's for $R^t(K, \mu)$ with $\lambda \rightarrow e_\lambda (f)$ analytic in $B(\lambda _0, \delta)$ for all $f \in R^t(K, \mu),$ then we say that $\lambda _0$ is an analytic bounded point evaluation (abpe) for $R^t(K, \mu).$ We denote the
set of abpe's for $R^t(K, \mu)$ by $abpe(R^t(K, \mu)).$ Clearly analytic bounded point evaluations are contained in the interior of $K.$

\cite{thomson} proves a remarkable structural theorem for $P^t(\mu ):$ 
There is a Borel partition $\{\Delta_i\}_{i=0}^\infty$ of $spt\mu$ such that the space $P^t(\mu |_{\Delta_i})$ contains no nontrivial characteristic functions and
 \[
 \ P^t(\mu ) = L^t(\mu |_{\Delta_0})\oplus \left \{ \oplus _{i = 1}^\infty P^t(\mu |_{\Delta_i}) \right \}.
 \]
Furthermore, if $U_i$ is the open set of analytic bounded point evaluations for
$P^t(\mu |_{\Delta_i})$ for $i \ge 1,$ then $U_i$ is a simply connected region and the closure of $U_i$ contains $\Delta_i.$

Because of Thomson's decomposition, the study of general $P^t(\mu )$ can be reduced to the case where $P^t(\mu )$ is irreducible (contains no nontrivial characteristic functions) and $abpe(P^t(\mu ))$ is a nonempty simply connected open set whose
closure contains $spt \mu.$ \cite{oy95} shows that one can use the Riemann Mapping Theorem to further reduce to the case where $abpe(P^t(\mu )) = \mathbb D.$ In this case, \cite{ars} obtained the following remarkable structural theorem.

\begin{ARSTheorem} \label{ARSTheorem}
Suppose that $\mu$ is supported in $\bar{\mathbb D}$ and is such that
$abpe (P^t(\mu )) = \mathbb D$ and $P^t(\mu )$ is irreducible, and that $\mu (\partial \mathbb D)> 0.$ Then:
\newline
a) If $f \in P^t(\mu )$ then the nontangential limit $f^*(z)$ of f exists for $\mu |_{\partial \mathbb D}$-
almost all $z,$ and $f^* = f |_{\partial \mathbb D}$ as elements of $L^t(\mu |_{\partial \mathbb D}).$
\newline
b) Every nonzero invariant subspace of $P^t(\mu )$ has index 1.
\end{ARSTheorem}

\cite{ce93} extends some results of Thomson's Theorem to the space $R^t(K,\mu ).$ \cite{b08} expresses $R^t(K,\mu )$ as a direct sum as the following:  With the assumption that the diameters of the components of $\mathbb C\setminus K$ are bounded away from zero, there exists a Borel partition $\{\Delta_i\}_{i=0}^\infty$ of $spt\mu$ and matching compact subsets $\{K_i\}_{i=0}^\infty$ of $K$ such that $\Delta_i \subset K_i$
and 
 \[
 \ R^t(K,\mu ) = L^t(\mu |_{\Delta_0})\oplus \left \{ \oplus _{i = 1}^\infty R^t(K_i, \mu |_{\Delta_i}) \right \}, \tag{1-1}
 \]
where for each $i\ge 1$ the corresponding summand $R^t(K_i, \mu |_{\Delta_i})$ is irreducible in the
sense that it contains no non-trivial characteristic function. Furthermore, if $U_i = abpe(R^t(K_i, \mu |_{\Delta_i}))$ for $i \ge 1,$ then $U_i$ is a connected region and the closure of $U_i$ contains $\Delta_i.$
The results includes both Thomson's theorem and results of \cite{ce93}.

It is evident that some restriction on the nature of $\mathbb C \setminus K$ is necessary in order ensure (1-1) to be
valid in general. Because of Brennan's decomposition under some additional conditions for $\mathbb C \setminus K,$ it is reasonable to assume, in the study of general $R^t(K, \mu ),$ that $R^t(K, \mu )$ is irreducible and $abpe(R^t(K,\mu ))$ is a nonempty connected open set whose closure contains $spt \mu.$ It is the purpose of this paper to  explore the boundary values of functions and indices of rationally invariant subspaces for $R^t(K, \mu )$ and to extend Aleman-Richter-Sundberg's Theorem.

Notice that it is possible for two compact sets, $K_1$ and $K_2,$ to contain the support of $\mu $ and satisfy $R^t(K_1, \mu ) = R^t(K_2, \mu ).$ Thus giving conditions on a compact set $K$ is inappropriate unless attention is focused on the smallest compact set which yields the same set of functions. Since $K \supset \sigma (S_\mu ),$ the spectrum of $S_\mu,$ $\sigma (S_\mu )$ is the smallest set. We will always assume that $K = \sigma (S_\mu ).$

For readability purpose, in section 2, we consider an important special case that the boundary of unbounded component of $\mathbb C\setminus K$ is the unit circle.  Proposition \ref{MProposition1}, which locally estimates the boundary values of Cauchy transform of an annihilating measure in the sense of capacitary density, plays a key role in proving Theorem \ref{MTheorem1} that extends Aleman-Richter-Sundberg's Theorem. As a consequence, our approach provides an alternative short proof of Aleman-Richter-Sundberg's Theorem. 
The main difficulty in their original proof, in \cite{ars}, is the proof of the following inequality:
 \[
 \ \underset{\lambda\rightarrow z}{\overline{\lim}} (1 - |\lambda |^2) ^{\frac{1}{t}} M_{\lambda } \le \dfrac{C}{h(z)^{\frac{1}{t}}} \tag{1-2}
 \]
nontangentially for $m$-almost all $z \in \partial \mathbb D,$ where $C$ is some constant. Our proof does not depend on the inequality (1-2). However, we will also develop a more general version of (1-2) in section 3 (see Theorem \ref{MTheorem4}). Proposition \ref{MProposition2}, which estimates the upper bound of Cauchy transform of an annihilating measure, is used to prove Theorem \ref{MTheorem2} that extends Theorem C in \cite{ars}.  

To facilitate the discussion of further results for more general $K,$ we provide the following example.
\begin{Example}
	Let $ 0 < \epsilon < \frac{1}{8},~ M =\{z: ~ -\frac{1}{2} < Re(z) < \frac{1}{2}, ~  Im(z) =0\},$ $U_n = \{z: ~ -\frac{1}{2} < Re(z) < \frac{1}{2}, ~ \frac{1}{2^n}(\frac{1}{2} - \epsilon )< Im(z) < \frac{1}{2^n}(\frac{1}{2} + \epsilon )\},$ and $L_n = \{z: ~ -\frac{1}{2} < Re(z) < \frac{1}{2}, ~ \frac{1}{2^n}(-\frac{1}{2} - \epsilon )< Im(z) < \frac{1}{2^n}(-\frac{1}{2} + \epsilon )\}.$ Let
	\[
	\ K_1 = \bar {\mathbb D} \setminus \left ( \cup _{n=1}^\infty U_n \right )
	\]
	and 
	\[
	\ K_2 = \bar {\mathbb D} \setminus \left ( (\cup _{n=1}^\infty U_n)\cup (\cup _{n=1}^\infty L_n) \right )
	\]
	Let $\mu$ and $\nu$ be positive finite measures with $spt(\mu )\subset K_1$ and $spt(\nu )\subset K_2$ so that $R^t(K_1,\mu)$ and $R^t(K_2,\nu)$ are irreducible. Suppose that $abpe(R^t(K_1,\mu)) = Int(K_1)$ (for example, $\mu = Area |_{Int(K_1)} + m|_{M},$ where $m|_{M}$ is Lebesgue measure on $M$) and $abpe(R^t(K_2,\nu)) = Int(K_2)$. By the Radon-Nikodym theorem, we can write $\mu = \mu _a + \mu _s $ and $\nu = \nu _a + \nu _s, $ where $\mu _a << m|_{M},$ $\mu _s \perp m|_{M},$ $\nu _a << m|_{M},$ and $\mu _s \perp m|_{M}.$    
\end{Example}
In this example, $M$ is the inner boundary of $K_i.$  It is natural to explore nontangential limits of functions of $R^t(K_1, \mu )$ on the inner boundary $M$ (from below) with respect to $\mu_a.$ What can we say about $R^t(K_2, \nu )?$ 

The purpose of section 3 is to investigate the boundary behaviors of the functions in $R^t(K, \mu )$ for the boundaries other than the unit circle in section 2. Theorem \ref{SBTheorem} proves if $R^t(K, \mu )$ is irreducible and there are 'big parts' of $\mathbb C \setminus K$ near 'both sides' of $E\subset \partial K,$ then $\mu (E) = 0.$ In the above example, the inner boundary $M$ of $K_2$ satisfies the property, so our result implies $\nu_a(E) = 0.$ Therefore, it is not needed to investigate the values of functions in $R^t(K_2, \nu )$ for the boundary $M.$ Theorem \ref{SBTheorem} can also be applied to those $K$ for which the diameters of all components of $\mathbb C\setminus K$ are bounded away from zero. For example, if $K$ in Theorem \ref{MTheorem1} or Theorem \ref{MTheorem2} satisfies the property, then the carrier of $\mu |_{\partial \mathbb D}$ is away from $\mathbb D\setminus K.$ In the case, the nontangential limits of functions in $R^t(K, \mu )$ can be defined with respect to $\mu |_{\partial \mathbb D}.$ Theorem \ref{MTheorem3} generalizes Theorem \ref{MTheorem1}. Finally, Theorem \ref{MTheorem4} generalizes the inequality (1-2) ((1.4) in \cite{ars}).

Before closing this section, we mention here a few related papers. For a compactly supported complex  measure $\nu$ of $\mathbb C,$ by estimating analytic capacity of the set $\{\lambda: |\mathcal C\nu (\lambda)| \ge c \},$ where $\mathcal C\nu$ is Cauchy transform of $\nu$ (see section 2 for definition),  \cite{b06}, \cite{ars}, and \cite{ARS10} provide interesting alternative proofs of Thomson's theorem. Both their proofs rely on X. Tolsa's deep results on analytic capacity. The author refines the estimations for Cauchy transform, in Lemma 4 of \cite{y17}, to study the bounded point evaluations for rationally multicyclic subnormal operators. Also see the work of \cite{a01}, \cite{a02}, \cite{ar97}, \cite{ms90}, \cite{msy99}, \cite{ot80}, \cite{ty95}, \cite{tr79a}, \cite{tr79b}, \cite{wy98}, \cite{y95a}, and \cite{y95b}.

\section{Outer boundary of $K$ is the unit circle}

In this section, we will concern the special cases where the outer boundary of $K$ is the unit circle $\partial \mathbb D.$ Consequently, we provide alternative proofs of Theorem A and Theorem C in \cite{ars}.

Let $\nu$ be a compactly supported finite measure on $\mathbb {C}.$ The Cauchy transform
of $\nu$ is defined by
\[
\ \mathcal C\nu (z) = \int \dfrac{1}{w - z} d\nu (w)
\]
for all $z\in\mathbb {C}$ for which
$\int \frac{d|\nu|(w)}{|w-z|} < \infty .$ A standard application of Fubini's
Theorem shows that $\mathcal C\nu \in L^s_{loc}(\mathbb {C} )$ for $ 0 < s < 2,$ in particular, it is
defined for area-almost all $z,$ and clearly $\mathcal C\nu$ is analytic in $\mathbb C_\infty \setminus spt \nu,$ where $\mathbb C_\infty = \mathbb C \cup \{\infty \}.$

For a compact $K \subset \mathbb C$ we
define the analytic capacity of $K$ by
\[
\ \gamma(K) = sup |f'(\infty)|
\]
where the sup is taken over those functions $f$ analytic in $\mathbb C_\infty \setminus K$ for which
$|f(z)| \le 1$ for all $z \in \mathbb C_\infty \setminus K,$ and
$f'(\infty) = \lim _{z \rightarrow \infty} z[f(z) - f(\infty)].$
The analytic capacity of a general $E \subset \mathbb C$ is defined to be 
\[
\ \gamma (E) = \sup \{\gamma (K) : K \subset E, ~K~ compact\}.
\]
Good sources for basic information about analytic
capacity are \cite{Gar72}, Chapter VIII of \cite{gamelin}, Chapter V of \cite{conway}, and \cite{tol14}.

A related capacity, $\gamma _+,$ is defined for $E \subset \mathbb C$ by
\[
\ \gamma_+(E) = sup \|\mu \|
\]
where the sup is taken over positive measures $\mu$ with compact support
contained in $E$ for which $\|\mathcal C\mu \|_{L^\infty (\mathbb C)} \le 1.$ Since $\mathcal C\mu$ is analytic in $\mathbb C_\infty \setminus spt \mu$ and $(\mathcal C \mu)'(\infty) = −\|\mu \|,$ we have
\[
\ \gamma _+(E) \le \gamma (E)
\]
for all $E \subset \mathbb C.$ \cite{Tol03} proves the astounding result (Tolsa's Theorem) that 
$\gamma_+$ and $\gamma$ are actually equivalent. 
 That is, there is an absolute constant $A_T$ such that
\[ 
\ \gamma (E) \le A_ T \gamma_+(E) \tag{2-1}
\]
for all $E \subset \mathbb C.$ The following semiadditivity of analytic capacity is a conclusion of Tolsa's Theorem.
\[
\ \gamma \left (\bigcup_{i = 1}^m E_i \right ) \le A_T \sum_{i=1}^m \gamma(E_i)\tag{2-2}
\]
where $E_1,E_2,...,E_m \subset \mathbb C.$

Let $\nu$ be a compactly supported finite measure on $\mathbb {C}.$ For $\epsilon > 0,$ $\mathcal C_\epsilon \nu$ is defined by
\[
\ \mathcal C_\epsilon \nu (z) = \int _{|w-z| > \epsilon}\dfrac{1}{w - z} d\nu (w),
\]
and the maximal Cauchy transform is defined by
 \[
 \ \mathcal C_* \nu (z) = \sup _{\epsilon > 0}| \mathcal C_\epsilon \nu (z) |.
 \]
The 1-dimensional radial maximal operator of $\nu$ (see also (2.7) in \cite{tol14}) is defined by
 \[
 \ M_R \nu (z) = \sup _{r > 0} \dfrac{|  \nu  | (B(z. r))}{r}.
 \]

\begin{Lemma}\label{TolsaLemma} 
There is an absolute positive constant $C_T,$ for $a > 0,$ we have
\newline
(1)  
\[
\ \gamma(\{\mathcal C_*\nu  \geq a\}) \le \dfrac{C_T}{a} \|\nu \|, \tag{2-3}
\]
\newline
(2) 
\[
\ m(\{M_R\nu  \geq a\}) \le \dfrac{C_T}{a} \|\nu \|. 
\]
In this case, if we define
 \[
 \ MV(\nu ) = \{e^{i\theta}:  M_R\nu (e^{i\theta}) = +\infty\}, \tag{2-4}
 \]
then $m(MV(\nu )) = 0.$ 
\end{Lemma}

{\bf Proof:} (1) follows from Proposition 2.1 of \cite{Tol02} and Tolsa's Theorem (2-1) (also see \cite{tol14} Proposition 4.16). Theorem 2.6 in \cite{tol14} implies (2). 
 
For $0 < \sigma < 1$ and $z \in \partial \mathbb D,$ we define the nontangential approach region $\Gamma _\sigma (z)$ to be the interior of the convex hull of $\{z\} \cup B(0,\sigma ).$ It is well known that the existence of nontangential limits on a set $E \subset \partial \mathbb D$ is independent of $\sigma$ up to sets of $m$-measure zero, so we will write $\Gamma (z) = \Gamma _{\frac{1}{2}}(z)$ a nontangential approach region. The following lemma is due to Lemma 1 in \cite{kt77}.

\begin{Lemma}\label{KTLemma}
Suppose $\nu$ is a finite positive measure supported on $\mathbb D,$ define 
 \[
 \ IV(\nu ) = \{e^{i\theta}:~ \underset{\Gamma (e^{i\theta})\ni \lambda \rightarrow e^{i\theta}}{\overline\lim}\int_{\mathbb D} \dfrac{1 - |\lambda|^2}{|1 - \bar \lambda z|^2} d \nu (z) > 0 \}\tag{2-5}
 \]
then $m(IV(\nu )) = 0.$
\end{Lemma}
 
For a finite compactly supported measure $\nu ,$ definite
 \[
 \ U(\nu ) = \{\lambda \in \mathbb C: ~ \int \dfrac{1}{|z - \lambda |} d|\nu |(z) < \infty\}.
 \]
Then $Area ((U(\nu ))^c) = 0.$

\begin{Lemma}\label{CauchyTLemma} 
Let $\nu$ be a finite  measure supported in $\bar {\mathbb D}$ and $| \nu | (\partial \mathbb D ) = 0.$ Let $1 < p \le \infty ,$ $q = \frac{p}{p-1},$ $f \in C (\bar {\mathbb D}),$ and $g \in L^{q} (| \nu |).$ Define
 \[
 \ EV(|g|^q|\nu | ) = MV (|g|^q|\nu | ) \cup IV (|g|^q|\nu | ) \tag{2-6}
 \]
where $MV (|g|^q|\nu | )$ and $IV (|g|^q|\nu | )$ are defined as in (2-4) and (2-5), respectively.
Suppose that $a > 0$ and $e^{i\theta} \in \partial \mathbb D \setminus EV(|g|^q|\nu | ),$ then there exist $\frac{3}{4} < r_\theta < 1,$ $E_\delta ^f \subset \bar B (e^{i\theta}, \delta ),$ and $\epsilon (\delta ) > 0,$ where $0 < \delta < 1 - r_\theta ,$ such that 
 \[
 \ \lim _{\delta \rightarrow 0} \epsilon (\delta ) = 0,
 \] 
 \[
 \ \gamma(E_\delta ^f) <\epsilon (\delta ) \delta ,
 \]
and for $|\lambda _0 -  e^{i\theta} | = \frac{\delta}{2} $ and $\lambda _0\in \Gamma (e^{i\theta }),$
\[
\ \left  |\mathcal C\left (\dfrac{(1 - \bar \lambda _0 z)^{\frac{2}{p}}}{(1 - |\lambda _0|^2)^{\frac{1}{p}}}fg\nu \right )(\lambda) - \mathcal C\left (\dfrac{(1 - \bar \lambda _0 z)^{\frac{2}{p}}}{(1 - |\lambda _0|^2)^{\frac{1}{p}}}fg\nu \right )(\frac{1}{\bar{\lambda} _0}) \right | \le a\|f\|_{L^{p} (| \nu |)}
\]
for all $\lambda\in (B (e^{i\theta}, \delta ) \setminus E_\delta ^f )\cap U(g\nu ) .$ Notice that $E_\delta ^f$ depends on $f$ and all other parameters are independent of $f.$ 
\end{Lemma}

{\bf Proof:} For $e^{i\theta} \in \partial \mathbb D\setminus EV(|g|^q|\nu | ),$ by Lemma \ref{TolsaLemma} and \ref{KTLemma}, we conclude that $m(EV(|g|^q|\nu | )) = 0,$ $M_1 = M_R (|g|^{q}|\nu |) (e^{i\theta}) < \infty ,$ and there exists $\frac{3}{4} < r_\theta < 1$ such that
 \[
 \  \left (\int_{\mathbb D}\dfrac{(1 - |\lambda _0|^2) |g|^q}{|1 - \bar \lambda _0 z|^2}d|\nu | \right )^{\frac{1}{q}} \le \frac{a}{256}\tag{2-7}
 \]
for $\delta < 1 - r_\theta.$ 
Let $\nu_\delta = \frac{\chi _{B (e^{i\theta}, N \delta )}}{(1 - \bar{\lambda} _0 z )^{1-\frac{2}{p}}\delta^{\frac{1}{p}}}fg \nu .$  For $\epsilon < \delta ,$ $N > 2 ,$ and $\lambda \in B (e^{i\theta} , \delta ),$ we get:
 \[
 \ 2(1 - |\lambda _0|) \le \delta \le 4(1 - |\lambda _0|),
 \]
 \[
 \ \bar B (\lambda , \epsilon) \subset B (e^{i\theta}, 2 \delta ) \subset B (e^{i\theta}, N\delta),
 \]
and
 \[ 
 \ \begin{aligned}
 \ & \left  |\mathcal C _\epsilon \left ((1 - \bar \lambda _0 z)^{\frac{2}{p}}\delta^{-\frac{1}{p}}fg\nu \right )(\lambda) - \mathcal C  \left ((1 - \bar \lambda _0 z)^{\frac{2}{p}}\delta^{-\frac{1}{p}}fg\nu \right )(\frac{1}{\bar{\lambda} _0}) \right | \\
 \ \le & \dfrac{|1 - \bar{\lambda} _0 \lambda |}{\delta^{\frac{1}{p}}} \left | \int _{|z - \lambda| > \epsilon}\dfrac{fgd\nu  }{(z - \lambda)(1 - \bar{\lambda}_0 z)^{1-\frac{2}{p}}} \right | + \left | \mathcal C  \left (\chi _{\bar B (\lambda, \epsilon)}\dfrac{(1 - \bar \lambda _0 z)^{\frac{2}{p}}}{\delta^{\frac{1}{p}}}fg\nu \right )(\frac{1}{\bar{\lambda} _0}) \right |\\
\ \le & 2\delta^{\frac{1}{q}}  \left | \int _{ B (e^{i\theta}, N\delta )^c}\dfrac{fgd\nu}{(z - \lambda)(1 - \bar{\lambda}_0 z)^{1-\frac{2}{p}}} \right | + 2\delta \left |\int _{|z - \lambda| > \epsilon}\dfrac{d\nu_\delta}{(z - \lambda)} \right | \\
 \ & + \int_{\bar B (\lambda, \epsilon)}\dfrac{\delta^{-\frac{1}{p}}}{|1 - \bar \lambda _0 z|^{1-\frac{2}{p}}}|fg|d|\nu | \\
 \ \le & 2\delta^{\frac{1}{q}} \sum_{k=0}^{\infty}\int _{2^kN\delta \le |z - e^{i\theta} | < 2^{k+1}N\delta} \dfrac{|f||g|d|\nu |}{|z - \lambda ||1 - \bar{\lambda} _0 z |^{1-\frac{2}{p}} }  + 2\delta |\mathcal C_\epsilon \nu _\delta (\lambda )| \\
 \ & + \int_{ B (e^{i\theta}, 2 \delta)}\dfrac{|1 - \bar \lambda _0 z|\delta^{-\frac{1}{p}}}{|1 - \bar \lambda _0 z|^{\frac{2}{q}}}|fg|d|\nu | \\
 \ \le & 2\delta^{\frac{1}{q}} \sum_{k=0}^{\infty} \dfrac{(2^{k+1}N\delta)^{\frac{1}{q}}(2^kN\delta  + 2\delta )^{\frac{2}{p}}}{(2^kN\delta  - \delta )(2^kN\delta  - 2\delta )} \left (\dfrac{\int _{B (e^{i\theta}, 2^{k+1}N\delta)}|g|^qd|\nu |}{2^{k+1}N\delta} \right )^{\frac{1}{q}} \|f\|_{L^{p} (| \nu |)} \\
\ &+ 2\delta \mathcal  C_* \nu _\delta (\lambda ) + 4\int_{ B (e^{i\theta}, 2 \delta)}\dfrac{\delta^{\frac{1}{q}}}{|1 - \bar \lambda _0 z|^{\frac{2}{q}}}|fg|d|\nu |\\
 \ \le & \dfrac{4(N+2)^{1+\frac{1}{p}}\sum_{k=0}^\infty 2^{\frac{-k}{q}}M_1^{\frac{1}{q}}}{(N-1)(N-2)} \|f\|_{L^{p} (| \nu |)} + 2\delta \mathcal C_* \nu _\delta (\lambda ) \\
\ & + 4 \|f\|_{L^{p} (| \nu |)} \left (\int_{\mathbb D}\dfrac{\delta |g|^q}{|1 - \bar \lambda _0 z|^2}d|\nu | \right )^{\frac{1}{q}}\\
\ \end{aligned}
 \]
Let 
 \[
 \ N = 6 + \left (\dfrac{256}{a} \sum_{k=0}^\infty 2^{\frac{-k}{q}}  \right)^q M_1, 
 \]
then together with (2-7), we get
 \[
 \ \begin{aligned}
 \ & \left  |\mathcal C _\epsilon \left ((1 - \bar \lambda _0 z)^{\frac{2}{p}}\delta^{-\frac{1}{p}}fg\nu \right )(\lambda) - \mathcal C  \left ((1 - \bar \lambda _0 z)^{\frac{2}{p}}\delta^{-\frac{1}{p}}fg\nu \right )(\frac{1}{\bar{\lambda} _0}) \right | \\
\ \le &\dfrac{a}{8} \|f\|_{L^{p} (| \nu |)} + 2\delta  \mathcal C_* \nu _\delta (\lambda )
\ \end{aligned} \tag{2-8}
 \]
for $\lambda \in B (e^{i\theta} , \delta ).$
Let
 \[
 \ E_\delta ^f = \{\lambda : \mathcal C_* \nu _\delta (\lambda ) \ge \frac{a\|f\|_{L^{p} (| \nu |)}}{16\delta } \} \cap \bar B (e^{i\theta}, \delta ),
 \]
From (2-3) and Holder's inequality, we get
 \[
 \ \begin{aligned}
 \ \gamma (E_\delta ^f) \le & \dfrac{16C_T\delta }{a\|f\|_{L^{p} (| \nu |)}} \int _{ B (e^{i\theta}, N\delta )}\dfrac{|f||g|d|\nu |}{|1 - \bar{\lambda} _0 z |^{1-\frac{2}{p}}\delta^{\frac{1}{p}}} \\
 \ \le & \dfrac{16C_T\delta }{a\|f\|_{L^{p} (| \nu |)}} \int _{ B (e^{i\theta}, N\delta )}\dfrac{|1 - \bar{\lambda} _0 z |\delta^{-\frac{1}{p}}|f||g|d|\nu |}{|1 - \bar{\lambda} _0 z |^{\frac{2}{q}}} \\
 \ \le & \dfrac{16C_T(N+2)\delta }{a\|f\|_{L^{p} (| \nu |)}} \int _{ B (e^{i\theta}, N\delta )}\dfrac{\delta^{\frac{1}{q}}|f||g|d|\nu |}{|1 - \bar{\lambda} _0 z |^{\frac{2}{q}}} \\
\ \le & \dfrac{64C_T(N+2) \delta}{a} \left(\int_{\mathbb D} \dfrac{1 - |\lambda _0 |^2}{|1 - \bar {\lambda} _0 z|^2} |g|^qd | \nu | \right ) ^ {\frac{1}{q}} 
 \ \end{aligned}\tag{2-9}
 \]
Set 
 \[
 \ \epsilon (\delta ) = \dfrac{65C_T(N+2)}{a} \left(\int_{\mathbb D} \dfrac{1 - |\lambda _0 |^2}{|1 - \bar {\lambda} _0 z|^2} |g|^qd | \nu | \right ) ^ {\frac{1}{q}}, 
 \]
then $\lim_{\delta \rightarrow 0 } \epsilon (\delta ) = 0$ and 
 \[
 \ \gamma (E_\delta ^f) <  \epsilon (\delta ) \delta.
 \]
From (2-8) and the definition of $E_\delta ^f,$ for $\lambda \in B (e^{i\theta}, \delta ) \setminus E_\delta ^f$ and $\epsilon < \delta ,$ we conclude that
 \[
 \ \begin{aligned}
 \ & \left  |\mathcal C _\epsilon \left (\dfrac{(1 - \bar \lambda _0 z)^{\frac{2}{p}}}{(1 - |\lambda _0|^2)^{\frac{1}{p}}}fg\nu \right )(\lambda) - \mathcal C\left (\dfrac{(1 - \bar \lambda _0 z)^{\frac{2}{p}}}{(1 - |\lambda _0|^2)^{\frac{1}{p}}}fg\nu \right )(\frac{1}{\bar{\lambda} _0}) \right | \\
 \ \le &4\left  |\mathcal C _\epsilon \left ((1 - \bar \lambda _0 z)^{\frac{2}{p}}\delta^{-\frac{1}{p}}fg\nu \right )(\lambda) - \mathcal C \left ((1 - \bar \lambda _0 z)^{\frac{2}{p}}\delta^{-\frac{1}{p}}fg\nu \right )(\frac{1}{\bar{\lambda} _0}) \right | \\
\  < & a \|f\|_{L^{p} (| \nu |)}
\ \end{aligned} 
\]
The lemma follows since the limit of $\mathcal C _\epsilon,$ when $\epsilon\rightarrow 0,$  exists for $ \lambda\in (B (e^{i\theta}, \delta ) \setminus E_\delta ^f )\cap U(g\nu ).$

\begin{Proposition}\label{MProposition1} Let $\nu$ be a finite complex  measure with support in $K \subset \bar {\mathbb D}.$  Suppose that $\nu \perp Rat(K)$ and $\nu | _{\partial \mathbb D} = hm$ ($m = \frac{d\theta}{2\pi}$).  Then for $b > 0$ and $m$-almost all $e^{i\theta}\in \partial \mathbb D,$  there exist $\frac{3}{4} < r_\theta < 1,$ $E_{\delta}\subset B(e^{i\theta}, \delta),$ and $\epsilon (\delta ) > 0,$ where $0 < \delta < 1 -r_\theta ,$  such that $\lim_{\delta\rightarrow 0}\epsilon (\delta ) = 0,$ $\gamma(E_\delta) < \epsilon (\delta ) \delta ,$ and 
\[
\ \left  |\mathcal C\nu (\lambda) - e^{-i\theta}h(e^{i\theta}) \right | \le b \tag{2-10}
\]
for all $\lambda\in (B (e^{i\theta}, \delta ) \cap \Gamma (e^{i\theta})\setminus E_\delta )\cap U(\nu ) .$
\end{Proposition}

{\bf Proof:} Let $\nu _1 = \nu | _{\mathbb D}$ and $\nu _2 = \nu | _{\partial \mathbb D} = hm.$ Using Plemelj's formula (see page 56 of \cite{cmr2006} or Theorem 8.8 in \cite{tol14}), we can find $E_1 \subset \partial \mathbb D$ with $m(E_1) = 0$ such that
 \[
 \ \lim_{ \Gamma (e^{i\theta})\ni z\rightarrow e^{i\theta}} \mathcal C \nu _2 (z) - \lim_{ \Gamma (e^{i\theta})\ni z\rightarrow e^{i\theta}} \mathcal  C \nu _2 (\frac{1}{\bar z}) = e^{-i\theta} h (e^{i\theta}) \tag{2-11}
 \]
for $e^{i\theta} \in \partial \mathbb D\setminus E_1.$ Set $E _0 = E_1\cup EV(|\nu_1|),$ where $EV(|\nu_1|)$ is defined as in (2-6) and $m(EV) = 0.$  

We now apply Lemma \ref{CauchyTLemma} for $p=\infty , ~ q = 1,~f = 1,~ g = 1,$ and $a=\frac{b}{2}.$ For $e^{i\theta} \in \partial \mathbb D\setminus E_0,$ there exist $\frac{3}{4} < r_\theta < 1,$ $E_{\delta}\subset B(e^{i\theta}, \delta),$ and $\epsilon (\delta ) > 0,$ where $0 < \delta < 1 -r_\theta ,$  such that $\lim_{\delta\rightarrow 0}\epsilon (\delta ) = 0,$ $\gamma(E_\delta) <(\delta ) \delta ,$ and for $\lambda _0 \in (\partial B (e^{i\theta}, \frac{\delta}{2} )) \cap \Gamma (e^{i\theta}),$
 \[
 \ \left | \mathcal C \nu _ 1(\lambda ) -  \mathcal  C \nu _ 1(\frac{1}{\bar \lambda_0})  \right | \le \frac{b}{2}
 \]
for all $\lambda\in (B (e^{i\theta}, \delta )\setminus E_\delta )\cap U(\nu ) .$ Moreover, from (2-11), $r_{\theta}$ can be chosen so that 
\[
 \ \left | \mathcal C \nu _2 (\lambda ) -  \mathcal  C \nu _2 (\frac{1}{\bar \lambda_0}) - e^{-i\theta} h (e^{i\theta}) \right | \le \frac{b}{2}
 \]
for $\lambda\in B (e^{i\theta}, \delta ) \cap \Gamma (e^{i\theta}).$ Since $\mathcal C\nu (\frac{1}{\bar{\lambda} _0}) = 0,$ we get
\[
\ \begin{aligned}
\ & \left  |\mathcal C\nu (\lambda) - e^{-i\theta}h(e^{i\theta}) \right | \\
\ \le & \left  |\mathcal C\nu _1 (\lambda) - \mathcal C\nu _1 (\frac{1}{\bar \lambda _0}) \right | + \left  |\mathcal C\nu_2 (\lambda) - \mathcal C\nu _2 (\frac{1}{\bar \lambda _0}) - e^{-i\theta}h(e^{i\theta}) \right | \\
\ \le b
\ \end{aligned}
\]
all $\lambda\in (B (e^{i\theta}, \delta ) \cap \Gamma (e^{i\theta})\setminus E_\delta )\cap U(\nu ) .$ The proposition is proved.

Let $R = \{ z: -1/2 < Re(z),Im(z) < 1/2 \}$ and $Q = \bar{\mathbb D}\setminus R.$ For a bounded Borel set
$E\subset \mathbb C$ and $1\le p \le \infty,$ $L^p(E)$ denotes the $L^p$ space with respect to the area measure $dA$ restricted to $E.$  The following Lemma is a simple application of Thomson's coloring scheme.

\begin{Lemma} \label{lemmaDSet}
There is an absolute constant $\epsilon _1 > 0$ with the
following property. If $\gamma (\mathbb D \setminus K) < \epsilon_1,$ then
\[
\ |f(\lambda ) | \le \|f\|_{L^\infty (Q\cap K)}
\]
for $\lambda \in R$ and $f \in A(\mathbb D),$ the uniform closure of $\mathcal P$ in $C(\bar {\mathbb D}).$
\end{Lemma}

{\bf Proof:} Let $S$ be a closed square whose edges are parallel to x-axis and y-axis. $S$ is defined to be light if $Area(S \cap K) = 0 .$ $S$ is heavy if it is not light. 

We now sketch our version of Thomson's coloring scheme for $Q$ with a given a positive integer $m.$ We refer the reader to \cite{thomson} and \cite{thomson3} section 2 for details.
 
For each integer $k > 3$ let $\{S_{kj}\}$ be an enumeration of the closed squares contained in $\mathbb C$ with edges of length $2^{-k}$
parallel to the coordinate axes, and corners at the points whose coordinates
are both integral multiples of $2^{-k}$ (except the starting square $S_{m1}$, see (3) below). 
In fact, Thomson's coloring scheme is just needed to be modified slightly as the following:

(1) Use our definition of a light $\epsilon$ square.

(2) A path to $\infty$ means a path to any point that is outside of $Q$ (replacing the polynomially convex hull of $\Phi$ by $Q$).

(3) The starting yellow square $S_{m1}$ in the $m$-th generation is $R.$ Notice that the length of $S_{m1}$ in $m$-th generation is $1$ (not $2^{-m}$).

We will borrow notations that are used in Thomson's coloring scheme such as $\{\gamma_n\}_{n\ge m}$ and $\{\Gamma_n\}_{n\ge m},$ etc. We denote 
 \[
 \ YellowBuffer_m = \sum _{k = m+1}^\infty k^2 2^{-k}.
 \]

Suppose the scheme terminates, in our setup, this means Thomson's coloring scheme reaches a square $S$ in $n$-th generation that is not contained in $Q.$ One can construct a polygonal path $P,$ which connects the centers of adjacent squares, from the center of a square (contained in $Q$) adjacent to $S$ to the center of a square adjacent to $R$ so that the orange (non green in Thomson's coloring scheme) part of length is no more than $YellowBuffer_m.$ Let $GP = \cup S_j,$ where $\{S_j\}$ are all light squares with  $P\cap S_j \ne \emptyset .$ By Tolsa's Theorem (2-2), we see
 \[
 \ \gamma (P) \le A_T (\gamma (Int(GP)) + YellowBuffer_m).
 \]
Since $P$ is a connected set, $\gamma (P) \ge \frac{0.1}{4}$ (Theorem 2.1 on page 199 of \cite{gamelin}).  We can choose $m $ to be large enough so that
 \[
 \ \gamma (GP) \ge \dfrac{1}{40A_T} - YellowBuffer_{m} = \epsilon_m > 0.
 \]
Now by Lemma 3 in \cite{b06} (or the proof of Case I of Lemma B in \cite{ars} on page 462-464), there exists a constant $\epsilon _0 > 0$ such that 
 \[
 \ \gamma (GP \setminus K) \ge \epsilon _0 \gamma (GP) \ge \epsilon _0 \epsilon _m.\tag{2-13}
 \]
So we have prove if the scheme terminates, then (2-13) holds.
 
Set $\epsilon_1 = \epsilon_0\epsilon_m.$ By assumption $\gamma (\mathbb D \setminus K) < \epsilon_1,$ we must have $\gamma (GP \setminus K) \le \gamma (\mathbb D \setminus K) < \epsilon_1. $ Therefore, the scheme will not terminate since (2-13) does not hold. In this case, one can construct a sequence of heavy barriers inside $Q,$ that is, $\{\gamma_n\}_{n\ge m}$ and $\{\Gamma_n\}_{n\ge m}$ are infinite.

Let $f\in A(\mathbb D),$ by the maximal modulus principle, we can find $z_n\in\gamma_n$ such that $|f(\lambda )| \le |f(z_n)|$ for $\lambda \in R.$ By the definition of $\gamma_n,$ we can find a heavy square $S_n$ with $z_n\in S_n\cap\gamma_n.$ Since $Area(S_n\cap K) > 0,$ we can choose $w_n\in S_n$ with $|f(w_n)| = \|f\|_{L^\infty (S_n\cap K)}.$ $\frac{f(w)-f(z_n)}{w-z_n}$ is analytic in $\mathbb D,$ therefore, by the maximal modulus principle again, we get
 \[
 \ \left | \dfrac{f(w_n)-f(z_n)}{w_n-z_n} \right | \le \sup_{w \in \gamma_{n+1}} \left | \dfrac{f(w)-f(z_n)}{w-z_n} \right | \le \dfrac{2\|f\|_{L^\infty (\mathbb D)}}{dist (z_n,\gamma_{n+1})}.
 \]
Therefore,
 \[
 \ |f(\lambda )| \le |f(z_n)| \le |f(w_n)| + \dfrac{2|z_n-w_n|\|f\|_{L^\infty (\mathbb D)}}{dist (z_n,\gamma_{n+1})} \le \|f\|_{L^\infty(Q\cap K)} + \dfrac{2\sqrt 2 2^{-n}\|f\|_{L^\infty (\mathbb D)}}{n^2 2^{-n}}
 \]
 for $\lambda \in R.$
The lemma follows by taking $n\rightarrow \infty .$

\begin{Corollary} \label{CorollaryDSet}
There is an absolute constant $\epsilon _1 > 0$ with the
following property. If $\lambda _0 \in \mathbb C, ~ \delta > 0,$ and $\gamma (B(\lambda _0 , \delta) \setminus K) < \epsilon_1\delta ,$ then
\[
\ |f(\lambda ) | \le \|f\|_{L^\infty (B(\lambda _0 , \delta)\cap K)}
\]
for $\lambda \in B(\lambda _0 , \frac{\delta}{2})$ and $f \in A(B(\lambda _0 , \delta)),$ the uniform closure of $\mathcal P$ in $C(\bar B(\lambda _0 , \delta)).$
\end{Corollary}

Now we assume that $R^t(K,\mu )$ is irreducible and $\Omega$ is a connected region satisfying:
 \[
 \ abpe(R^t(K,\mu )) = \Omega ,~ K = \bar \Omega, ~ \Omega \subset \mathbb D,~ \partial \mathbb D \subset \partial\Omega . \tag{2-14}
 \]
It is well known that, in this case, $\mu |_{\partial \mathbb D} << m. $
 So we assume $\mu |_{\partial \mathbb D} = hm.$ 

For $0 < \delta < 1$ and $e^{i\theta}\in \partial \mathbb D,$ define $\Gamma ^\delta _\sigma (e^{i\theta}) = \Gamma _\sigma (e^{i\theta})\cap B(e^{i\theta}, \delta).$ In order to define a nontangential limit of a function in $R^t(K,\mu )$ at $e^{i\theta} \in \partial \Omega,$ one needs $\Gamma ^\delta _\sigma (e^{i\theta}) \subset \Omega$ for some $\delta.$ Therefore, we define the strong outer boundary of $\Omega$ as the following:
 \[
 \ \partial _{so, \sigma} \Omega = \{e^{i\theta} \in \partial \Omega: ~\exists 0<\delta<1,~ \Gamma _{\sigma}^\delta (e^{i\theta}) \subset \Omega \}, ~ \partial _{so} \Omega = \partial _{so, \frac{1}{2}} \Omega. \tag{2-15}
 \]
It is known that $\partial _{so, \sigma} \Omega$ is a Borel set (i.e., see Lemma 4 in \cite{ot80}) and $m(\partial _{so, \sigma_1} \Omega \setminus \partial _{so, \sigma _2} \Omega ) = 0$ for $\sigma _1 \ne \sigma _2.$ From Theorem \ref{SBTheorem} in section 3, if $R^t(K,\mu )$ is irreducible and the diameters of all components of $\mathbb C\setminus K$ are bounded away from zero, then $\mu (\partial \mathbb D\setminus \partial _{so} \Omega ) = 0.$ This means that the carrier of $\mu | _{\partial \mathbb D} $ is a subset of $\partial _{so} \Omega $ and the nontangential limit of a function at $e^{i\theta}\in \partial \mathbb D\setminus \partial _{so} \Omega $ is not defined.

 From Lemma VII.1.7 in \cite{conway}, we find a function $G \in R^t(K,\mu )^\perp \subset L^{t'}(\mu )$ such that $G(z) \ne 0$ for $\mu$-almost every $z.$ Every $f \in R^t(K,\mu )$ is analytic on $\Omega$ and 
\[
\ f(\lambda ) \mathcal C(G\mu ) (\lambda ) = \int \dfrac{f(z)}{z - \lambda }G(z) d\mu(z) = \mathcal C(fG\mu ) (\lambda )\tag{2-16}
\]
for $\lambda \in \Omega \cap U(G\mu ).$     

\begin{Theorem}\label{MTheorem1}
Suppose that $\mu$ is a finite positive measure supported in $K$ and is such that
$abpe(R^t(K,\mu )) = \Omega $ and $R^t(K,\mu )$ is irreducible, where $\Omega$ is a connected region satisfying (2-14), $\mu |_{\partial \mathbb D} = hm,$ and $\mu (\partial _{so} \Omega ) > 0.$ Then:
\newline
(a) If $f \in R^t(K,\mu )$ then the nontangential limit $f^*(z)$ of $f$ exists for $\mu |_{\partial _{so} \Omega}$-
almost all z, and $f^* = f |_{\partial _{so} \Omega}$ as elements of $L^t(\mu |_{\partial _{so} \Omega}).$
\newline
(b) Every nonzero rationally invariant subspace $M$ of $R^t(K,\mu )$ has index 1, that is, $dim(M / (S_\mu - \lambda _0) M) = 1,$ for $\lambda _0\in \Omega.$
\newline
If the diameters of all components of $\mathbb C \setminus K$ are bounded away from zero, then by Theorem \ref{SBTheorem} (in section 3), the above $\partial _{so} \Omega$ can be replaced by $\partial \mathbb D.$ 
\end{Theorem}

{\bf Proof:} (a) Let $1 > \epsilon > 0$ and $\epsilon _0 = \frac{\epsilon _1}{32A_T},$ where $\epsilon _1$ is as in Lemma \ref{lemmaDSet} and $A_T$ is from (2-2). For $f \in R^t(K,\mu ),$ from Proposition \ref{MProposition1}, we see that for $\mu$-almost all $e^{i\theta}\in \partial _{so} \Omega$ with $\Gamma ^{r_0}(e^{i\theta}) \subset \Omega$ and $G(e^{i\theta})h(e^{i\theta} \ne 0,$ $b = \frac{|G(e^{i\theta})h(e^{i\theta})|}{2(1+|f(e^{i\theta})|)} \epsilon > 0,$ there exist $\max (r_0, \frac{3}{4} ) < r_\theta < 1,$ $E_{\delta}^1\subset B(e^{i\theta}, \delta),$ $E_{\delta}^2\subset B(e^{i\theta}, \delta),$ and $\epsilon (\delta ) > 0,$ where $0 < \delta < 1 -r_\theta ,$  such that $\lim_{\delta\rightarrow 0}\epsilon (\delta ) = 0,$ $\gamma(E_\delta ^1) < \epsilon (\delta ) \delta ,$ $\gamma(E_\delta ^2) < \epsilon (\delta ) \delta ,$ 
\[
\ \left  |\mathcal C(G\mu ) (\lambda) - G(e^{i\theta})e^{-i\theta}h(e^{i\theta}) \right | \le b
\]
for all $\lambda\in (B (e^{i\theta}, \delta ) \cap \Gamma (e^{i\theta})\setminus E_\delta ^1)\cap U(G\mu ),$ and
\[
\ \left  |\mathcal C(fG\mu ) (\lambda) - f(e^{i\theta}) G(e^{i\theta}) e^{-i\theta}h(e^{i\theta}) \right | \le b
\]
for $\lambda\in (B (e^{i\theta}, \delta ) \cap \Gamma (e^{i\theta})\setminus E_\delta ^2)\cap U(G\mu ).$ Now choose $\delta$ small enough so that $\epsilon(\delta ) < \epsilon_0.$ Set $E_\delta  = E_\delta ^1 \cup E_\delta ^2,$ then from the semi-additivity (2-2), we get
 \[
 \ \gamma (E_\delta ) \le A_T (\gamma (E_\delta ^1) + \gamma (E_\delta ^2)) < \epsilon _1 \dfrac{\delta}{16} .
 \]
Therefore, by (2-16), for $\lambda\in (B (e^{i\theta}, \delta ) \cap \Gamma (e^{i\theta})\setminus E_\delta )\cap U(G\mu ),$
 \[ 
 \ \begin{aligned}
 \ & |f (\lambda) - f(e^{i\theta})| \\
\ \le &\left | \dfrac{\mathcal C(fG\mu ) (\lambda) - f(e^{i\theta})\mathcal C(G\mu ) (\lambda)}{\mathcal C(G\mu ) (\lambda)}\right |  \\
 \ \le & \dfrac{2|\mathcal C(fG\mu ) (\lambda) - f(e^{i\theta}) G(e^{i\theta}) e^{-i\theta}h(e^{i\theta})|}{|G(e^{i\theta})h(e^{i\theta})|} + \dfrac{2|\mathcal C(G\mu ) (\lambda) -  G(e^{i\theta}) e^{-i\theta}h(e^{i\theta})| |f(e^{i\theta})|}{|G(e^{i\theta})h(e^{i\theta})|}\\
 \ \le &  \epsilon .
 \ \end{aligned}
 \]
For $\lambda _0 \in (\partial B (e^{i\theta}, \frac{\delta}{2} )) \cap \Gamma _{\frac{1}{4}}(e^{i\theta}),$ we see that $B (\lambda _0, \frac{\delta}{16}) \subset B (e^{i\theta}, \delta ) \cap \Gamma (e^{i\theta}).$ 
 Using Lemma \ref{lemmaDSet} for $f  - f(e^{i\theta}),$ we get 
 \[
 \ |f (\lambda) - f(e^{i\theta})| \le \|f  - f(e^{i\theta}) \| _{L^\infty (B (\lambda _0, \frac{\delta}{16}) \setminus E_\delta)} \le \epsilon 
 \]
for every $\lambda \in B (\lambda_0, \frac{\delta }{32}).$ Hence,
\[
 \ \lim_{ \Gamma _{\frac{1}{4}}(e^{i\theta})\ni\lambda\rightarrow e^{i\theta}}f(\lambda ) = f(e^{i\theta} ).
 \]

We turn to prove (b). Let $M$ be a nonzero rationally invariant subspace of $R^t(K,\mu ).$ Without loss of generality, we assume $\lambda_0 = 0$ and $0\in \Omega.$ We must show that $dim(M/S_\mu M) = 1.$ Let $n$ be the smallest integer such that $f(z) = z^n f_0(z)$ for every $f\in M$ and there exists $g\in M$ with $g(z) = z^n g_0(z)$ and $g_0(0) 
\ne 0.$ We only need to show $\frac{f(z) - \frac{f_0(0)}{g_0(0)} g(z)}{z}\in M.$ To do this, it is suffice to show that for $\phi\in M^\perp \subset L^{t'} (\mu)$, the function
 \[
 \ \Phi (\lambda ) = \int \dfrac{g(\lambda )f(z) - f(\lambda )g(z)}{z - \lambda } \phi (z) d\mu (z),
 \]
which is analytic in $\Omega,$ is identically zero. In fact, the proof is similar to that of (a). Let $E \subset \partial _{so} \Omega$ so that for $e^{i\theta} \in E,$ $f$ and $g$ have nontangential limits at $e^{i\theta},$ and $h(e^{i\theta}) > 0. $ By Theorem \ref{MTheorem1} (a), $m(E) > 0.$  For $1 > \epsilon > 0$ and $\epsilon _0 = \frac{\epsilon _1}{32A_T},$ applying Proposition \ref{MProposition1} for $f\phi\mu , g\phi\mu$ since $f\phi\mu , g\phi\mu \perp Rat(K)$ and Theorem \ref{MTheorem1} (a) for $f$ and $g,$  we see that for $e^{i\theta} \in E$ with $\Gamma ^{r_0} (e^{i\theta}) \subset \Omega$ and  $b = \frac{1}{(1 + |f(e^{i\theta})| + |g(e^{i\theta})| ) (1 + |\phi(e^{i\theta})| h(e^{i\theta}))} \epsilon,$ there exist $\max (r_0, \frac{3}{4} ) < r_\theta < 1,$ $E_{\delta}^1\subset B(e^{i\theta}, \delta),$ $E_{\delta}^2\subset B(e^{i\theta}, \delta),$ and $\epsilon (\delta ) > 0,$ where $0 < \delta < 1 -r_\theta ,$  such that $\lim_{\delta\rightarrow 0}\epsilon (\delta ) = 0,$ $\gamma(E_\delta ^1) < \epsilon (\delta ) \delta ,$ $\gamma(E_\delta ^2) < \epsilon (\delta ) \delta ,$ 
\[
\ \left  |\mathcal C(f\phi\mu ) (\lambda) - f(e^{i\theta})\phi(e^{i\theta})e^{-i\theta}h(e^{i\theta}) \right | \le b
\]
for all $\lambda\in (B (e^{i\theta}, \delta ) \cap \Gamma (e^{i\theta})\setminus E_\delta ^1)\cap U(G\mu ),$
\[
\ \left  |\mathcal C(g\phi\mu ) (\lambda) - g(e^{i\theta})\phi(e^{i\theta})e^{-i\theta}h(e^{i\theta}) \right | \le b,
\]
for all $\lambda\in (B (e^{i\theta}, \delta ) \cap \Gamma (e^{i\theta})\setminus E_\delta ^2)\cap U(G\mu ),$ 
$|f(\lambda) - f(e^{i\theta})| < b$ and $|g (\lambda) - g(e^{i\theta})| < b$ on $B (e^{i\theta}, \delta ) \cap \Gamma (e^{i\theta}).$
Choose $\delta$ small enough so that $\epsilon (\delta ) < \epsilon_0.$ Set $E_\delta  = E_\delta ^1 \cup E_\delta ^2,$ then by the semi-additivity (2-2) again, we have $\gamma (E_\delta ) < \epsilon _1 \frac{\delta}{16} .$ 
Therefore, for all $\lambda\in (B (e^{i\theta}, \delta ) \cap \Gamma (e^{i\theta})\setminus E_\delta )\cap U(G\mu ),$
 \[ 
 \ \begin{aligned}
 \ & |\Phi (\lambda) | \\
\ \le & |g(\lambda)| \left  |\mathcal C(f\phi\mu ) (\lambda) - f(e^{i\theta})\phi(e^{i\theta})e^{-i\theta}h(e^{i\theta}) \right | + |f(\lambda)|   |\mathcal C(g\phi\mu ) (\lambda) \\ 
\ &- g(e^{i\theta})\phi(e^{i\theta})e^{-i\theta}h(e^{i\theta}) | + |f(\lambda)g(e^{i\theta}) - g(\lambda)f(e^{i\theta})| |\phi(e^{i\theta})| h(e^{i\theta})\\
 \ \le & (b + |f(e^{i\theta})| + |g(e^{i\theta})| ) (1 + |\phi(e^{i\theta})| h(e^{i\theta})) b\\
 \ \le &  \epsilon .
 \ \end{aligned}
 \]
For $\lambda _0 \in (\partial B (e^{i\theta}, \frac{\delta}{2} )) \cap \Gamma _{\frac{1}{4}}(e^{i\theta}),$ $B (\lambda _0, \frac{\delta}{16}) \subset B (e^{i\theta}, \delta ) \cap \Gamma (e^{i\theta}).$
Using Lemma \ref{lemmaDSet} for $\Phi,$ we get 
 \[
 \ |\Phi (\lambda )| \le \|\Phi \| _{L^\infty (B (\lambda _0, \frac{\delta}{16}) \setminus E_\delta)} \le \epsilon 
 \]
for every $\lambda \in B (\lambda_0, \frac{\delta }{32}).$ Hence,
\[
 \ \lim_{ \Gamma _{\frac{1}{4}}(e^{i\theta})\ni\lambda\rightarrow e^{i\theta}} \Phi (\lambda ) = 0.
 \]
Let $V = \cup _{e^{i\theta}\in E}\Gamma _{\frac{1}{4}}^\delta (e^{i\theta}).$ Since $m(E) > 0,$ there exists a connected component $V_0$ of $V$ with $m(\partial V_0 \cap \partial \mathbb D) > 0.$  $\partial V_0$ is a rectifiable Jordan curve and $\Phi (\lambda )$ is analytic in $V_0.$ Therefore $\Phi (\lambda ) = 0$ since $\Omega$ is a connected region. This completes the proof.

\begin{Proposition}\label{MProposition2} Let $\mu$ be a finite positive  measure with support in $ K \subset\bar {\mathbb D}$  and $\mu | _{\partial \mathbb D} = hm.$  Let $1 < p <\infty, ~ q = \frac{p}{p-1}, ~ f\in C(\bar{\mathbb D}), ~ g \in L^q (\mu ),$ and $fg\mu \perp Rat(K).$ Then for $0 < \beta < \frac{1}{16},$ $b > 0,$ and $m$-almost all $e^{i\theta}\in \partial \mathbb D,$ there exist $\frac{3}{4} < r_\theta < 1,$ $E_{\delta}^f \subset B(e^{i\theta}, \delta),$ and $\epsilon (\delta ) > 0,$ where $0 < \delta < 1 -r_\theta ,$ such that $\lim_{\delta\rightarrow 0}\epsilon (\delta ) = 0,$ $\gamma(E_\delta ^f) < \epsilon (\delta ) \delta ,$ and for $\lambda _0 \in (\partial B (e^{i\theta}, \frac{\delta}{2} )) \cap \Gamma _{\frac{1}{4}}(e^{i\theta}),$
\[
\ \left  |\mathcal C \left (\dfrac{(1 - \bar \lambda _0 z)^{\frac{2}{p}}}{ (1-|\lambda _0|^2)^{\frac{1}{p}}}fg\mu \right )(\lambda) \right | \le \left(b + \dfrac{1+4\beta}{1-4\beta} \left ( \int _{\partial \mathbb D} \dfrac{1 - |\lambda _0|^2}{|1 - \bar \lambda _0 z |^2} |g|^qd\mu \right )^{\frac{1}{q}} \right ) \|f\|_{L^p(\mu )} \tag{2-17}
\]
for all $\lambda\in (B (\lambda _0, \beta \delta ) \setminus E_\delta^f )\cap U(g\mu ) .$
\end{Proposition}

{\bf Proof:} Let $\nu = \mu | _{\mathbb D}.$ We now apply Lemma \ref{CauchyTLemma} for $p, ~ q,~f,~ g, $ and $a=b.$ For $e^{i\theta} \in \partial \mathbb D\setminus EV(|g|^q\nu)$ (as in (2-6) and $m(EV(|g|^q\nu)) = 0$), there exist $\frac{3}{4} < r_\theta < 1,$ $E_{\delta}^f \subset B(e^{i\theta}, \delta),$ and $\epsilon (\delta ) > 0,$ where $0 < \delta < 1 -r_\theta ,$ such that $\lim_{\delta\rightarrow 0}\epsilon (\delta ) = 0,$ $\gamma(E_\delta ^f) < \epsilon (\delta ) \delta ,$ and for $\lambda _0 \in (\partial B (e^{i\theta}, \frac{\delta}{2} )) \cap \Gamma _{\frac{1}{4}}(e^{i\theta}),$
\[
\ \left | \mathcal C \left (\dfrac{(1 - \bar \lambda _0 z)^{\frac{2}{p}}}{ (1-|\lambda _0|^2)^{\frac{1}{p}}}fg\nu \right )(\lambda ) - \mathcal C \left (\dfrac{(1 - \bar \lambda _0 z)^{\frac{2}{p}}}{ (1-|\lambda _0|^2)^{\frac{1}{p}}}fg\nu \right )(\frac{1}{\bar{\lambda} _0})\right | \le b  \|f\|_{L^p(\mu )}
 \]
for all $\lambda\in (B (e^{i\theta }, \delta ) \setminus E_\delta^f )\cap U(g\mu ) .$
  \[
\ \mathcal C \left (\dfrac{(1 - \bar \lambda _0 z)^{\frac{2}{p}}}{ (1-|\lambda _0|^2)^{\frac{1}{p}}}fg\mu \right )(\frac{1}{\bar{\lambda} _0}) = 0
\]
since $fg\mu \perp Rat(K).$ From Lemma \ref{CauchyTLemma}, for all $\lambda\in (B (\lambda _0, \beta \delta ) \setminus E_\delta^f )\cap U(g\mu ),$ we get
 \[
 \ \begin{aligned}
 \ & \left  |\mathcal C \left (\dfrac{(1 - \bar \lambda _0 z)^{\frac{2}{p}}}{ (1-|\lambda _0|^2)^{\frac{1}{p}}}fg\mu \right )(\lambda) \right | \\
 \ \le & \left  |\mathcal C \left (\dfrac{(1 - \bar \lambda _0 z)^{\frac{2}{p}}}{ (1-|\lambda _0|^2)^{\frac{1}{p}}}fg\nu \right )(\lambda) - \mathcal C \left (\dfrac{(1 - \bar \lambda _0 z)^{\frac{2}{p}}}{ (1-|\lambda _0|^2)^{\frac{1}{p}}}fg\nu \right )(\frac{1}{\bar{\lambda} _0})\right | \\
 \ & + \left  |\int_{\partial \mathbb D} \left ( \dfrac{1}{z - \lambda} - \dfrac{1}{z - \frac{1}{\bar{\lambda} _0}} \right )\dfrac{(1 - \bar \lambda _0 z)^{\frac{2}{p}}}{ (1-|\lambda _0|^2)^{\frac{1}{p}}}fg\mu\right | \\
 \ \le & b \|f\|_{L^p(\mu )} + \int_{\partial \mathbb D} \dfrac{|1 - \lambda \bar{\lambda} _0| }{|z - \lambda |} \dfrac{ (1 - |\lambda _0|^2)^{-\frac{1}{p}} }{ |1 - \bar \lambda _0 z|^{1 - \frac{2}{p}}} |fg| d \mu \\
 \ \le & b \|f\|_{L^p(\mu )} + \dfrac{1+4\beta}{1-4\beta} \int_{\partial \mathbb D} \dfrac{ (1 - |\lambda _0|^2)^{\frac{1}{q}}}{ |1 - \bar \lambda _0 z|^{\frac{2}{q}}} |fg| d \mu 
 \ \end{aligned}
 \]
 where the last step follows from
 \[
 \ \dfrac{|1 - \lambda \bar{\lambda} _0| }{|z - \lambda |} \le \dfrac{1 - |\lambda _0|^2 + |\lambda _0||\lambda - \lambda _0| }{|z - \lambda _0 | - |\lambda - \lambda _0|} \le \dfrac{(1+4\beta)(1 - |\lambda _0|^2)}{|z - \lambda _0 | - 4\beta (1 - |\lambda _0|)} \le \dfrac{(1+4\beta)(1 - |\lambda _0|^2)}{(1 - 4\beta )|1 - \bar\lambda _0 z|}
 \]
for $z\in\partial \mathbb D.$ The corollary now follows from Holder's inequality.

\begin{Theorem}\label{MTheorem2} Suppose that $\mu$ is a finite positive measure supported in $K$ and is such that
$abpe(R^t(K,\mu )) = \Omega $ and $R^t(K,\mu )$ is irreducible, where $\Omega$ is a connected region satisfying (2-14), $\mu |_{\partial \mathbb D} = hm,$ and $\mu (\partial _{so} \Omega ) > 0.$
 Then for $t > 1,$ 
 \[
 \ \lim_{ \Gamma _{\frac{1}{4}}(e^{i\theta})\ni\lambda \rightarrow e^{i\theta}} (1 - |\lambda |^2)^{\frac{1}{t}} M_\lambda = \dfrac{1}{h(e^{i\theta})^{\frac{1}{t}}}
 \]
for $\mu$-almost all $e^{i\theta}\in \partial _{so} \Omega.$ If the diameters of all components of $\mathbb C \setminus K$ are bounded away from zero, then by Theorem \ref{SBTheorem} (in section 3), the above $\partial _{so} \Omega$ can be replaced by $\partial \mathbb D.$
\end{Theorem}

{\bf Proof:} By Proposition \ref{MProposition1} and \ref{MProposition2}, for $\mu$-almost all $e^{i\theta}\in \partial _{so} \Omega$ with $G(e^{i\theta})h(e^{i\theta}) \ne 0$ and $\Gamma ^{r_0}(e^{i\theta})  \subset \Omega ,$  $0 < \beta < \frac{1}{16},$  $b > 0,$ and $f\in Rat(K),$ there exist $max(\frac{3}{4},  r_0) < r_\theta < 1,$ $E_{\delta} \subset B(e^{i\theta}, \delta),$ $E_{\delta}^f \subset B(e^{i\theta}, \delta),$ and $\epsilon (\delta ) > 0,$ where $0 < \delta < 1 -r_\theta ,$ such that $\lim_{\delta\rightarrow 0}\epsilon (\delta ) = 0,$ $\gamma(E_\delta ) < \epsilon (\delta ) \delta ,$ $\gamma(E_\delta ^f) < \epsilon (\delta ) \delta ,$ 
\[
\ \left  |\mathcal C(G\mu ) (\lambda) - e^{-i\theta}G(e^{i\theta})h(e^{i\theta}) \right | \le b \tag{2-18}
\]
for all $\lambda\in (B (e^{i\theta}, \delta ) \cap \Gamma (e^{i\theta})\setminus E_\delta )\cap U(G\mu ) .$
and 
\[
\ \left  |\mathcal C \left (\dfrac{(1 - \bar \lambda _0 z)^{\frac{2}{t}}}{ (1-|\lambda _0|^2)^{\frac{1}{t}}}fG\mu \right )(\lambda) \right | \le \left(b + \dfrac{1+4\beta}{1-4\beta} \left ( \int _{\partial \mathbb D} \dfrac{1 - |\lambda _0|^2}{|1 - \bar \lambda _0 z |^2} |G|^{t'}d\mu \right )^{\frac{1}{t'}} \right ) \|f\|_{L^t(\mu )} \tag{2-19}
\]
for $\lambda _0 \in \partial B (e^{i\theta}, \frac{\delta}{2} ) \cap \Gamma _{\frac{1}{4}}(e^{i\theta})$ and all $\lambda\in (B (\lambda _0, \beta \delta ) \setminus E_\delta^f )\cap U(G\mu ) .$ From semi-additivity of (2-2), we get
 \[
 \ \gamma (E_\delta \cup E_\delta^f) \le A_T(\gamma (E_\delta ) + \gamma ( E_\delta^f)) \le 2A_T \epsilon (\delta ) \delta .
 \] 
Let $\delta $ be small enough so that $\epsilon (\delta ) < \frac{\beta}{2A_T}\epsilon _1,$ where $\epsilon _1$ is as in Corollary \ref{CorollaryDSet}. From (2-16), (2-18), and (2-19), for $\lambda _0 \in \partial B (e^{i\theta}, \frac{\delta }{2}) \cap \Gamma _{\frac{1}{4}}(e^{i\theta})$ and all $\lambda\in (B (\lambda _0, \beta \delta ) \setminus (E_{\delta } \cup E_{\delta }^f) )\cap U(G\mu ) ,$ we have the following calculation:
 \[
 \ | 1 - \bar \lambda _0 \lambda | \ge 1- |\bar \lambda _0|^2 - |\lambda - \lambda _0||\bar \lambda _0| \ge 1- |\bar \lambda _0|^2 - \beta \delta |\lambda _0|
 \]
and 
 \[
 \ \begin{aligned}
 \ (1-|\lambda _0|^2)^{\frac{1}{t}} |f(\lambda ) | \le & \dfrac{| (1 - \bar \lambda _0 \lambda )^{\frac{2}{t}}(1-|\lambda _0|^2)^{-\frac{1}{t}}f(\lambda ) |}{(1  - \beta \frac{\delta  |\lambda _0|}{1-|\lambda _0|^2})^{\frac{2}{t}}} \\
 \ = & \dfrac{1}{(1  - \beta \frac{\delta  |\lambda _0|}{1-|\lambda _0|^2})^{\frac{2}{t}}} \left |\dfrac{  \mathcal C \left (\dfrac{(1 - \bar \lambda _0 z)^{\frac{2}{t}}}{ (1-|\lambda _0|^2)^{\frac{1}{t}}}fG\mu \right )(\lambda) }{\mathcal C(G\mu ) (\lambda)} \right | \\
\ \le & \dfrac{ b + \frac{1+4\beta}{1-4\beta} \left ( \int _{\partial \mathbb D} \frac{1 - |\lambda _0|^2}{|1 - \bar \lambda _0 z |^2} |G|^{t'}d\mu \right )^{\frac{1}{t'}} }{(1-4\beta)^{\frac{2}{t}}(|G(e^{i\theta})|h(e^{i\theta}) - b)} \|f\|_{L^t(\mu )}.
\ \end{aligned}
\] 
Since $\gamma (E_{\delta }\cup E_{\delta }^f) < \epsilon _1 \delta ,$ from Corollary \ref{CorollaryDSet}, we conclude
 \[
 \ M_{\lambda _0} \le \sup _{\underset{\|f\|_{L^t(\mu )} = 1}{f\in Rat(K)}}|f(\lambda _0)| \le \sup _{\underset{\|f\|_{L^t(\mu )} = 1}{f\in Rat(K)}}\|f\|_{L^\infty (B (\lambda _0, \beta \delta ) \setminus (E_{\delta } \cup E_{\delta }^f) }  
 \]
for $\lambda _0 \in \partial B (e^{i\theta}, \frac{\delta}{2} ) \cap \Gamma _{\frac{1}{4}}(e^{i\theta}).$  Hence,
 \[
 \ \underset{\Gamma _{\frac{1}{4}}(e^{i\theta})\ni \lambda _0 \rightarrow e^{i\theta}}{\overline\lim} (1-|\lambda _0|^2)^{\frac{1}{t}} M_{\lambda _0} \le \dfrac{ b + \frac{1+4\beta}{1-4\beta} |G(e^{i\theta})|(h(e^{i\theta}))^{\frac{1}{t'}} }{(1-4\beta)^{\frac{2}{t}}(|G(e^{i\theta})|h(e^{i\theta}) - b)}
 \]
since $\frac{1 - |\lambda _0|^2}{|1 - \bar \lambda _0 z |^2}$ is the Poisson kernel. Taking $b\rightarrow 0$ and $\beta\rightarrow 0, $ we get
 \[
 \ \lim_{\Gamma _{\frac{1}{4}}(e^{i\theta}) \ni \lambda \rightarrow e^{i\theta}} (1 - |\lambda |^2)^{\frac{1}{t}} M_\lambda \le \dfrac{1}{h(e^{i\theta})^{\frac{1}{t}}}.
 \]
The reverse inequality is from \cite{kt77} (applying Lemma \ref{KTLemma} to testing function $(1 - \bar \lambda _0 z )^{-\frac{2}{t}}$). This completes the proof.
	
\section{Boundary values of $R^t(K,\mu)$ for certain $K$}
	
	In this section, we are concerning the boundary behaviors of functions in $R^t(K,\mu)$ near the boundary of $K$ (not necessarily outer boundary as in last section), in particular, the inner boundary of $K.$ Our approach in estimating Cauchy transform, in section 2, is concentrating on the local behavior of the transform. This makes it possible to extend our methodology to more general $K.$ In order to apply our approach, the following requirements are needed.
	\newline
	(A) Plemelj's formula must hold for the boundary points under consideration;
	\newline
	(B) Lemma \ref{TolsaLemma} (2) and Lemma \ref{KTLemma} shall be extended.
	
	For (A), it is known that Plemelj's formula holds for a Lipschitz graph (see Theorem 8.8 in \cite{tol14}). So we will restrict our attention to the boundary of $K$ which is a part of a Lipschitz graph although Plemelj's formula may hold for more general rectifiable curves.

	We define the open cone (with vertical axis) 
	 \[
	 \ \Gamma (\lambda, \alpha ) = \{z \in \mathbb C :~ |Re(z) -  Re(\lambda )| < \alpha |Im(z) -  Im(\lambda )| \}, 
	 \]
	and the half open cones 
	 \[
	 \ \Gamma ^+ (\lambda, \alpha ) = \{z \in \Gamma (\lambda, \alpha ) :~ Im(z) >  Im(\lambda )\}
 \]
	and
	 \[
	 \ \Gamma ^- (\lambda, \alpha ) = \{z \in \Gamma (\lambda, \alpha ) :~ Im(z) <  Im(\lambda ) \}.
	 \]
	Set $\Gamma _\delta ^+ (\lambda, \alpha ) = B(\lambda ,  \delta) \cap \Gamma ^+ (\lambda, \alpha )$ and $\Gamma _\delta ^- (\lambda, \alpha ) = B(\lambda ,  \delta) \cap \Gamma ^- (\lambda, \alpha ).$ $\Gamma ^+ (\lambda, \alpha )$ (or $\Gamma _\delta ^+ (\lambda, \alpha )$) is called upper cone. $\Gamma ^- (\lambda, \alpha )$ (or $\Gamma _\delta ^- (\lambda, \alpha )$) is called lower cone.
		
	Let $A : ~ \mathbb R \rightarrow \mathbb R$ be a Lipschitz function and let $LG$ be its graph. Observe that if $\alpha < \frac{1}{\|A'\|_\infty},$ then, for every $\lambda\in LG,$
	$\Gamma ^+ (\lambda, \alpha )\subset \{z\in \mathbb C:~ Im(z) > A(Re(z))\}$ and $\Gamma ^- (\lambda, \alpha )\subset \{z\in \mathbb C:~ Im(z) < A(Re(z))\}.$
On the graph Γ of $A,$ we consider the usual complex measure
	 \[
	 \ dz_{LG} = \dfrac{1 + iA'(Re(z))}{ (1 + A'(Re(z))^2)^{\frac{1}{2}}} d\mathcal H^1 |_{LG} = (L(z))^{-1}d\mathcal H^1 |_{LG}\tag{3-1}
	 \]
	where $\mathcal H^1$ is one dimensional Hausdorff measure. Notice that $|L(z)| = 1.$ For $1 \le p < \infty $ and $f \in L^p(\mathcal H^1 |_{LG}),$ the nontangential limits 
	\[
	\ \mathcal C_+ (f dz_{LG}(\lambda) =  \lim_{\Gamma ^+ (\lambda, \alpha ) \ni z\rightarrow \lambda} \mathcal C (f dz_{LG}(z)	
	\]
	and 
	\[
	\ \mathcal C_{-} (f dz_{LG}(\lambda) =  \lim_{\Gamma ^{-} (\lambda, \alpha ) \ni z\rightarrow \lambda} \mathcal C (f dz_{LG}(z)	
	\]
	exist $\mathcal H^1 |_{LG}$-almost everywhere. Moreover,
	 \[
	 \ \dfrac{1}{2\pi i} \mathcal C_+ (f dz_{LG}(\lambda) - \dfrac{1}{2\pi i}C_- (f dz_{LG}(\lambda) = f(\lambda ) \tag{3-2}
	 \]
(see Theorem 8.8 in \cite{tol14}).

Suppose that $R^t(K,\mu )$ is irreducible and $\Omega$ is a connected region satisfying:
 \[
 \ abpe(R^t(K,\mu )) = \Omega ,~ K = \bar \Omega . \tag{3-3}
 \]	
Let $G \in R^t(K,\mu )^\perp \subset L^{t'}(\mu )$ such that $G(z) \ne 0$ for $\mu$-almost every $z.$

In order to apply our approach, we need to impose some constraints on $K$ and define type I and II boundaries for $K.$
Upper cone $\Gamma ^+(\lambda, \alpha )$ (or lower cone $\Gamma ^-(\lambda, \alpha )$) is outer for $\lambda \in LG\cap \partial K$ if there exist $\delta _\lambda, ~ \epsilon _\lambda > 0$ such that for every $\delta < \delta _\lambda,$
 \[
 \ B(\lambda _\delta, \epsilon _\lambda \delta)  \subset K^c \cap \Gamma _\delta ^+(\lambda, \alpha ) ~(\text{or } B(\lambda _\delta, \epsilon _\lambda \delta)  \subset K^c\cap \Gamma _\delta ^-(\lambda, \alpha )). \tag{3-4}
 \]
 $\lambda \in LG\cap \partial K$ is a type I boundary point of $LG\cap \partial K$ if either upper cone $\Gamma ^+(\lambda, \alpha )$ or lower cone $\Gamma ^-(\lambda, \alpha )$ is outer. The type I boundary $\partial _{I,\alpha}^{LG} K$ is the set of all type I boundary points of $LG\cap \partial K.$ For example, if $V$ is a component of $K$ and $\partial V$ is a Lipschitz graph, then
$\partial V$ is a type I boundary. 

Upper cone $\Gamma ^+(\lambda, \alpha )$ (or lower cone $\Gamma ^-(\lambda, \alpha )$) is inner for $\lambda \in LG\cap \partial K$ if there exists $\delta > 0$ such that 
\[
 \ \Gamma _\delta ^+(\lambda, \alpha ) \subset \Omega ~(\text{or } \Gamma _\delta ^-(\lambda, \alpha ) \subset \Omega ).
 \]
$\lambda \in LG\cap \partial K$ is a type II boundary point of $LG\cap \partial K$ if $\lambda$ is type I and either upper cone $\Gamma ^+(\lambda, \alpha )$ or lower cone $\Gamma ^-(\lambda, \alpha )$ is inner.
The type II boundary $\partial _{II,\alpha}^{LG} K$ is the set of all type II boundary points of $LG\cap \partial K.$ The strong outer boundary of $\Omega$ defined in the section 2 is type II boundary of $K.$

Without loss of generality, for type I boundary point $\lambda,$ we usually assume upper cone $\Gamma ^+(\lambda, \alpha )$ is outer, and for type II boundary point $\lambda,$ we usually assume lower cone $\Gamma ^-(\lambda, \alpha )$ is inner. 
 
\begin{Lemma}
Both $\partial _{I,\alpha}^{LG} K$ and $\partial _{II,\alpha}^{LG} K$ are Borel sets.
\end{Lemma} 

{\bf Proof:} Let $\epsilon _0 = \frac{1}{n}$ and define $A_{nm}$ to be the set of $\lambda \in LG\cap \partial K$ such that for every $0 < \delta < \frac{1}{m},$ there exists $\lambda _0$ with 
 \[
 \ B(\lambda _0,\epsilon _0\delta ) \subset K^c \cap \Gamma _\delta ^+(\lambda, \alpha). 
 \] 
One sees that $A_{nm}$ is a closed set and $\partial _{I,\alpha}^{LG} K = \cup A_{nm}.$ If we define $B_{nmk}$ to be the set of $\lambda \in A_{nm}$ such that $\Gamma _{\frac{1}{k}} ^- (\lambda, \alpha ) \subset Int(K),$ then it is straightforward to verify that $B_{nmk}$ is a closed set and $\partial _{II,\alpha}^{LG} K = \cup B_{nmk}.$ 

It is easy to verify  that $\mathcal H^1 |_{LG} (\partial _{I,\alpha_1}^{LG} K \setminus \partial _{I,\alpha_2}^{LG} K) = 0$ and $\mathcal H^1 |_{LG} (\partial _{II,\alpha_1}^{LG} K \setminus \partial _{II,\alpha_2}^{LG} K) = 0$ for $\alpha_1 \ne \alpha_2.$ Therefore, we will fix $0 < \alpha < \frac{1}{\|A'\|_\infty}$ and use $\partial _{I}^{LG} K,~ \partial _{II}^{LG} K$ for $\partial _{I,\alpha}^{LG} K,~ \partial _{II,\alpha}^{LG} K,$ respectively.

For (B), Lemma \ref{lemmaBasic} and Corollary \ref{CorollaryKTGen} below extend Lemma \ref{TolsaLemma} (2) and Lemma \ref{KTLemma}. From now on, we use $LG$ for a fixed Lipschitz graph as above.

\begin{Lemma}\label{lemmaBasic}
Let $\nu$ be a finite complex measure with compact support. Suppose $\nu$ is singular to $\mathcal H^1 |_{LG}$ ($|\nu |\perp \mathcal H^1 |_{LG}$). Then
\newline
(1)
 \[ 
 \ \mathcal H^1 |_{LG} (\{\lambda  : M_R\nu (\lambda )  \geq a\}) \le \dfrac{C}{a} \|\nu \|
 \]
where $C$ is an absolute constant. In this case,
 \[
 \ \mathcal H^1 |_{LG} (\{\lambda  : M_R\nu (\lambda )  = \infty\}) = 0 \tag{3-5}.
 \]
\newline
(2)
 \[
 \ \mathcal H^1 |_{LG} (\{\lambda  :\underset{\delta\rightarrow 0}{\overline\lim} \dfrac{|\nu |(B(\lambda , \delta )}{\delta} > 0 \}) = 0. \tag{3-6}
 \]

\end{Lemma}

{\bf Proof:} As the same as Lemma \ref{TolsaLemma} (2), (1) follows from Theorem 2.6 in \cite{tol14}. 

(2) Let $E_0$ be a Borel set such that $\mathcal H^1 |_{LG} (E_0) = 0$ and $|\nu | (E_0 ^c) = 0$ (since $|\nu |\perp \mathcal H^1 |_{LG}$). Let $\epsilon ,~ \eta > 0$ and let $E \subset \{\lambda  :\underset{\delta\rightarrow 0}{\overline\lim} \dfrac{|\nu |(B(\lambda , \delta )}{\delta} > \frac{1}{N} \} \cap E_0^c$ be a compact subset. Let $O$ be an open set containing $E$ with $|\nu | (O) < \eta .$ Let $x\in E,$ then there exists $ 0 < \delta _x < \frac{\epsilon}{3}$ such that $|\nu |(B(x, \delta _x )) \ge \frac{1}{N} \delta _x$ and $B(x, \delta _x ) \subset O.$ Since $E \subset \cup_{x\in E} B(x, \delta _x),$ we can choose a finite subset $\{x_i\}_{i=1}^n $ so that $E\subset \cup_{i=1}^n B(x_i, \delta _{x_i}).$ From $3r$-covering theorem (see Theorem 2.1 in \cite{tol14}), we can further select a subset $\{x_{i_j}\}_{j=1}^m$ such that $\{B(x_{i_j}, \delta _{x_{i_j}})\}$ are disjoint and 
 \[
 \ E \subset \cup_{i=1}^n B(x_i, \delta _{x_i}) \subset \cup _{j=1}^m  B(x_{i_j}, 3\delta _{x_{i_j}}).
 \] 
Therefore,
 \[
 \ \mathcal H^1_\epsilon (E) \le 3\sum _{j=1}^m \delta _{x_{i_j}} \le \dfrac{ 3}{N} \sum _{j=1}^m |\nu |(B(x_{i_j}, \delta _{x_{i_j}})) = \dfrac{ 3}{N}  |\nu |(\cup _{j=1}^mB(x_{i_j}, \delta _{x_{i_j}})) \le \dfrac{ 3}{N}  |\nu | (O) < \dfrac{ 3}{N} \eta . 
 \] 
This implies $\mathcal H^1 |_{LG}(E) = 0.$ The lemma is proved.

\begin{Corollary}\label{CorollaryKTGen}
Let $\nu$ be a positive finite compactly supported measure on $\mathbb C$ and $\nu$ is singular to $\mathcal H^1 |_{LG}$ ($\nu \perp \mathcal H^1 |_{LG}$). For $\mathcal H^1 |_{LG}$-almost all $w\in LG,$ if there exists $\delta_w, ~\epsilon _w > 0$ such that 
\[
 \ B(\lambda _\delta, \epsilon _w \delta ) \subset (spt(\nu ))^c\cap B(w, \delta )
 \]
for $0 < \delta < \delta_w,$ then 
 \[
 \ \lim _{\delta \rightarrow 0} \int \dfrac{\delta}{|z - \lambda _\delta|^2} d \nu (z) = 0.\tag{3-7}
 \]
\end{Corollary}

{\bf Proof:} From (3-5) and (3-6), we assume that
 \[
 \ M_R(w) < \infty,~ \underset{\delta\rightarrow 0}{\lim} \dfrac{\nu (B(w , \delta )}{\delta} = 0. 
 \]
Hence, for $N > 2,$
 \[
 \ \begin{aligned}
 \ & \int \dfrac{\delta}{|z - \lambda _\delta|^2} d \nu (z) \\
\ \le & \int _{B(w,N\delta )} \dfrac{\delta}{|z - \lambda _\delta|^2} d \nu (z) + \int _{B(w,N\delta )^c} \dfrac{\delta}{|z - \lambda _\delta|^2} d \nu (z) \\
 \ \le & \dfrac{N}{\epsilon _w^2} \dfrac{\nu (B(w , N\delta )}{N\delta} + \sum _{k = 0}^\infty \int _{2^kN\delta \le |z-w| < 2^{k+1}N\delta } \dfrac{\delta}{|z - \lambda _\delta|^2} d \nu (z) \\
 \ \le & \dfrac{N}{\epsilon _w^2} \dfrac{\nu (B(w , N\delta )}{N\delta} + \sum _{k = 0}^\infty \dfrac{2^{k+1}N\delta ^ 2}{(2^kN\delta - \delta)^2} \dfrac{\nu (B(w , 2^{k+1}N\delta )}{2^{k+1}N\delta} \\
 \ \le & \dfrac{N}{\epsilon _w^2} \dfrac{\nu (B(w , N\delta )}{N\delta} + \dfrac{4N}{(N-1)^2}M_R(w) 
 \ \end{aligned}
 \]
The second term is small for $N$ large and for a given $N,$ the first term is small if $\delta$ is small enough. Therefore, (3-7) holds. 

Now we state our generalized version of Lemma \ref{CauchyTLemma} below. Notice that there is no corresponding function $(1-\bar\lambda _0 z)^\frac{2}{p}$ for a boundary point $w$ of an arbitrary $K,$ in particular, for an inner boundary point $w$.

\begin{Lemma}\label{CauchyTLemmaGen} 
Let $\nu$ be a finite  measure supported in $K$ and $| \nu | \perp \mathcal H^1 |_{\partial _I^{LG} K}.$ Let $1 < p \le \infty ,$ $q = \frac{p}{p-1},$ $f \in C (K),$ and $g \in L^{q} (| \nu |).$ Define
 \[
 \ EVG(|g|^q|\nu | ) = \{\lambda \in \partial _I^{LG} K: ~ M_R(|g|^q||\nu |)(\lambda ) = \infty\text{ or }\underset{\delta \rightarrow 0} {\overline\lim}\int \dfrac{\delta |g|^q}{|z - \lambda _\delta|^2} d |\nu |(z) > 0 \}
 \]
where $\lambda _\delta$ is defined as in (3-4). Then $\mathcal H^1 |_{\partial _I^{LG} K}(EVG(|g|^q|\nu | ))=0$ (Lemma \ref{lemmaBasic} and
Corollary \ref{CorollaryKTGen}). Suppose that $a > 0,$ $w \in \partial _I^{LG} K \setminus EV(|g|^q|\nu | ),$ and upper cone $\Gamma _\delta ^+(w, \alpha )$ is outer, then there exist $\delta_w >0,$ $E_\delta ^f \subset \bar B (w, \delta ),$ and $\epsilon (\delta ) > 0,$ where $0 < \delta < \delta_w ,$ such that 
$\lim _{\delta \rightarrow 0} \epsilon (\delta ) = 0,$ $\gamma(E_\delta ^f) <\epsilon (\delta ) \delta ,
$
and 
\[
\ \left  |\mathcal C\left (fg\nu \right )(\lambda) - \mathcal C\left (fg\nu \right )(\lambda _\delta) \right | \le a \delta ^{-\frac{1}{p}}\|f\|_{L^{p} (| \nu |)}
\]
for all $\lambda\in (B (w, \delta ) \setminus E_\delta ^f )\cap U(g\nu  ) .$ Notice that $E_\delta ^f$ depends on $f$ and all other parameters are independent of $f.$ 
\end{Lemma}

{\bf Proof:} We just need to make the following slight modifications to the proof of Lemma \ref{CauchyTLemma}:
 \newline 
(1) Replace $\frac{1}{\bar \lambda _0}$ by $\lambda _\delta .$ 
\newline 
(2) Use Lemma \ref{lemmaBasic} (1) instead of Lemma \ref{TolsaLemma} (2) and use Corollary \ref{CorollaryKTGen} instead of  Lemma \ref{KTLemma}. 
\newline 
(3) Replace  $\nu_\delta$ by  $\nu_\delta = \frac{\delta^{\frac{1}{p}}\chi _ {B(w, N\delta )}}{z - \lambda _\delta}fg\nu .$ 
\newline
(4) (2-8) becomes
 \[
 \ \delta^{\frac{1}{p}}|\mathcal C_\epsilon \nu (\lambda ) - \mathcal C \nu (\lambda _\delta ) | \le \dfrac{a}{2} \|f\|_{L^{p} (| \nu |)} + 2\delta C_*\nu _\delta (\lambda ).
 \]
\newline
(5) Define
 \[
 \ E_\delta ^f = \{\lambda : C_*\nu _\delta (\lambda ) > \dfrac{a\|f\|_{L^{p} (| \nu |)}}{4\delta } \} \cap \bar B(w, \delta ).
 \]
(2-9) becomes
 \[
 \ \gamma (E_\delta ) \le \dfrac{4C_T\delta }{a\|f\|_{L^{p} (| \nu |)}} \int _{B(w,N \delta )} \dfrac{\delta^{\frac{1}{p}}|fg|d|\nu |}{|z - \lambda _\delta |} < \epsilon(\delta) \delta .
 \]   
where $\epsilon_w$ is as in (3-4) and
 \[
 \ \epsilon(\delta) = \dfrac{5(N+1)C_T}{a\epsilon_w}\left(  \int \dfrac{\delta |g|^qd|\nu |}{|z - \lambda _\delta |^2} \right )^\frac{1}{q}.
 \]	
The proof is completed.

\begin{Proposition}\label{MProposition1Gen} Let $\nu$ be a finite complex  measure with support in $K.$  Suppose that $\nu \perp Rat(K)$ and $\nu = \nu_a + \nu_s$ is the Radon-Nikodym decomposition with respect to $\mathcal H^1 |_{\partial _I^{LG} K},$ where $\nu_a = \frac{1}{2\pi} h\mathcal H^1 |_{\partial _I^{LG} K}$ and $\nu_s \perp \mathcal H^1 |_{\partial _I^{LG} K}.$ Suppose upper cone $\Gamma _\delta ^+(w, \alpha )$ is outer for $w\in \partial _I^{LG} K.$
Then for $b > 0$ and $\mathcal H^1 |_{\partial _I^{LG} K}$-almost all $w\in \partial _I^{LG} K,$  there exist $\delta _w > 0,$ $E_{\delta}\subset B(w, \delta),$ and $\epsilon (\delta ) > 0,$ where $0 < \delta < \delta _w,$  such that $\lim_{\delta\rightarrow 0}\epsilon (\delta ) = 0,$ $\gamma(E_\delta) < \epsilon (\delta ) \delta ,$ and 
\[
\ \left  |\mathcal C\nu (\lambda) - L(w) h(w) \right | \le b
\]
for all $\lambda\in (\Gamma _\delta ^-(w, \alpha )\setminus E_\delta )\cap U(\nu ) .$
\end{Proposition}

{\bf Proof:} We just need to replace Plemelj's formula (2-11) in the proof of Proposition \ref{MProposition1} by (3-2).

The following Lemma is from Lemma B in \cite{ars} (also see Lemma 3 in \cite{y17}).

\begin{Lemma} \label{lemmaARS}
There are absolute constants $\epsilon _1 > 0$ and $C_1 < \infty$ with the
following property. For $R > 0,$ let $E \subset \bar B(\lambda _0, R)$ with 
$\gamma (E) < R\epsilon_1.$ Then
\[
\ |p(\lambda)| \le \dfrac{C_1}{R^2} \int _{\bar B(\lambda _0, R)\setminus E} |p| \frac{dA}{\pi}
\]
for all $\lambda\in B(\lambda _0, \frac{R}{2})$ and $p \in A(B(\lambda _0, R)),$ the uniform closure of $\mathcal P$ in $C(\bar B(\lambda _0, R)).$
\end{Lemma} 

Set
\[
 \ a(\alpha) = \frac{1}{8}\sin(\frac{\tan^{-1}(\alpha)}{2}). \tag{3-8}
 \] 
Clearly, for $\lambda _0 \in\Gamma _\delta ^-(w, \frac{\alpha }{2}) \cap (\partial B(w, \frac{\delta}{2})),$
 \[
 \ B(\lambda _0, 2a(\alpha)\delta) \subset \Gamma _\delta ^-(w, \alpha ).
 \]
The following theorem indicates that the carrier of $\mu _a,$ for irreducible $R^t(K,\mu ),$ does not intersect the boundary points for which both upper and lower cones contain a big portion of $\mathbb C\setminus K.$

\begin{Theorem}\label{SBTheorem}
Suppose that $\mu$ is a finite positive measure supported in $K$ and is such that $abpe(R^t(K,\mu )) = \Omega$ and $R^t(K,\mu )$ is irreducible, where $\Omega$ satisfies (3-3). Suppose that upper cone $\Gamma _\delta ^+(w, \alpha )$ is outer for all $w\in \partial _I^{LG} K$ and $\mu = \mu_a + \mu_s$ is the Radon-Nikodym decomposition with respect to $\mathcal H^1 |_{\partial _I^{LG} K},$ where $\mu_a = \frac{1}{2\pi} h\mathcal H^1 |_{\partial _I^{LG} K}$ and $\mu_s \perp \mathcal H^1 |_{\partial _I^{LG} K}.$
\newline
(a) Define
 \[
 \ E = \{w\in \partial _I^{LG} K : \underset{\delta\rightarrow 0}{\overline{\lim}}\dfrac{\gamma (\Gamma _\delta^- (w,\alpha ) \setminus K)}{\delta} >0 \}, 
 \]
then $\mu _a (E) = 0.$
\newline
(b) If the diameters of all components of $\mathbb C \setminus K$ are bounded away from zero, then 
 \[
 \ \mu _a (\partial _I^{LG} K \setminus \partial _{II}^{LG} K ) = 0. 
 \]
\end{Theorem}

{\bf Proof:} (a) Let $G \in R^t(K,\mu )^\perp$ and $G(z) \ne 0~ \mu~ a.e.$ as above. Suppose $\mu _a (E) > 0, $ then there exists $w\in E$ such that 
 
(1) $G(w)h(w) \ne 0.$ 

(2) Proposition  \ref{MProposition1Gen} holds for $w,$ that is, for $b = \frac{|G(w)|h(w)}{2},$ there exist $\delta _w > 0,$ $E_{\delta}\subset B(w, \delta),$ and $\epsilon (\delta ) > 0,$ where $0 < \delta < \delta _w,$  such that $\lim_{\delta\rightarrow 0}\epsilon (\delta ) = 0,$ $\gamma(E_\delta) < \epsilon (\delta ) \delta ,$ and 
\[
\ \left  |\mathcal C(G\mu ) (\lambda) - L(w) G(w) h(w) \right | \le b \tag{3-9}
\]
for all $\lambda\in (\Gamma _\delta ^-(w, \alpha )\setminus E_\delta )\cap U(G\mu ) .$

(3) There are a sequence of $\{\delta_n\}$ with $\delta _n \rightarrow 0$ and $\epsilon _0 > 0$ such that
 \[
 \ \gamma (\Gamma _{\delta _n}^- (w,\alpha ) \setminus K) \ge \epsilon _0 \delta _n.
 \]
Choose $N$ large enough so that $\epsilon (\delta _N ) < \frac{\epsilon _0}{2}.$ For $\lambda \in \Gamma _{\delta _N}^- (w,\alpha ) \setminus K,$ we see that $\lambda \in U(G\mu )$ and (3-9) does not hold since $\mathcal C(G\mu) (\lambda) = 0.$ That implies
 \[
 \ \Gamma _{\delta _N}^- (w,\alpha )  \setminus K \subset E_{\delta _N}.
 \] 
Hence,
 \[
 \ \gamma (\Gamma _{\delta _N}^- (w,\alpha ) \setminus K) \le \gamma (E_{\delta _N}) \le \dfrac{\epsilon_0}{2}\delta _N. 
 \]
This contradicts (3).

We now turn to prove (b). Let $lb > 0$ be less than the diameters of all components of $\mathbb C \setminus K.$ 
Let $E_1$ be the set of $w\in \partial _I^{LG} K$ such that there exists a sequence of $\{\delta_n\}$ with $\delta_n\rightarrow 0$ and $\Gamma _{\delta _n}^- (w,\frac{\alpha}{2} ) \cap \partial K \ne \emptyset .$ For a given $w\in E_1,$ there exists a component $V_n$ of $\mathbb C\setminus K$ so that $\Gamma _{\delta _n}^- (w,\frac{\alpha}{2} ) \cap V_n \ne \emptyset .$
Let $\lambda _n \in \Gamma _{\delta _n}^- (w,\frac{\alpha}{2} ) \cap V_n,$ then
 \[
 \ B(\lambda _ n, a(\alpha )\delta _n) \cap V_n \subset \Gamma _{2\delta _n}^- (w,\alpha ) \setminus K,
 \]
where $a(\alpha )$ is defined as in (3-8).
Hence, 
 \[
 \ \begin{aligned}
 \ \dfrac{1}{4}\min (a(\alpha ) \delta _n, lb) \le &\dfrac{1}{4} diameter (B(\lambda _ n, a(\alpha )\delta _n) \cap V_n) \\
 \ \le &\gamma(B(\lambda _ n, a(\alpha )\delta _n) \cap V_n) \\
 \le &\gamma(\Gamma _{2\delta _n}^- (w,\alpha ) \setminus K),
 \end{aligned}
 \]  
where the second inequality is implied by Theorem 2.1 on page 199 of \cite{gamelin}. This implies
 \[
 \ \underset{r\rightarrow 0}{\overline{\lim}}\dfrac{\gamma (\Gamma _{\delta}^- (w,\alpha ) \setminus K)}{\delta} \ge \dfrac{a(\alpha )}{8}.
 \] 
So $E_1\subset E,$ from (a), we conclude $\mu _a (E_1) = 0.$ We have shown that $\Gamma _{\delta}^- (w,\frac{\alpha}{2} ) \cap \partial K = \emptyset $  for $w\in \partial _I^{LG} K\setminus E_1$ with $\mu _a (E_1) = 0$ as $\delta$ is close to zero enough, in this case, $\Gamma _{\delta}^- (w,\frac{\alpha}{2} ) \subset Int(K) .$ 

Let $w\in \partial _I^{LG} K\setminus E_1$ so that we can apply Proposition  \ref{MProposition1Gen} for $w$ and $b = \frac{ |G(w)h(w)|}{2}.$ There exist $\delta _w > 0,$ $E_{\delta}\subset B(w, \delta),$ and $\epsilon (\delta ) > 0,$ where $0 < \delta < \delta _w,$  such that $\lim_{\delta\rightarrow 0}\epsilon (\delta ) = 0,$ $\gamma(E_\delta) < \epsilon (\delta ) \delta ,$ and
\[
\ \left  |\mathcal C(G\mu ) (\lambda) \right | \ge \dfrac{ |G(w)h(w)|}{2} \tag{3-10}
\]
for all $\lambda\in (\Gamma _\delta ^-(w, \alpha )\setminus E_\delta )\cap U(G\mu ) .$ Now choose $\delta$ to be small enough so that $\epsilon (\delta ) < a(\frac{\alpha}{2})\epsilon _1,$ where $\epsilon _1$ is as in Lemma \ref{lemmaARS} and $a(\alpha)$ is defined in (2-8). Let $\lambda _0 \in \Gamma _{\delta} ^-(w, \frac{\alpha}{4} )$ with $|\lambda _0 - w| = \frac{\delta}{2},$ then $B(\lambda _0, a(\frac{\alpha}{2})\delta) \subset \Gamma _{\delta} ^-(w, \frac{\alpha}{2} ) \subset  Int(K),$ where $\delta$ is small enough. Since $\gamma (B(\lambda _0, a(\frac{\alpha}{2})\delta) \cap E_\delta ) < \epsilon_1 a(\frac{\alpha}{2})\delta , $ from Lemma \ref{lemmaARS} and (3-10). we conclude $\lambda \in B(\lambda _0, \frac{a(\frac{\alpha}{2})\delta}{2}),$
 \[
 \begin{aligned}
 \ |r(\lambda)| \le & \dfrac{C_1}{\pi (a(\frac{\alpha}{2})\delta)^2} \int _{B(\lambda _0, a(\frac{\alpha}{2})\delta) \setminus E_\delta } |r(z)| dA(z) \\
 \ \le & \dfrac{2C_1}{ \pi  |G(w)h(w)|a(\frac{\alpha}{2})^2 \delta ^2} \int _{B(\lambda _0, a(\frac{\alpha}{2})\delta) \setminus E_\delta } |\mathcal C(rG\mu ) (z)| dA(z) \\
\ \le & \dfrac{C_1}{ \pi  |G(w)h(w)| a(\frac{\alpha}{2})^2 \delta ^2} \int \int _{B(\lambda _0, a(\frac{\alpha}{2})\delta)} \dfrac{1}{|z-\lambda|} dA(z)|rG|d\mu  (\lambda)  \\
 \le & \dfrac{C_2}{\delta} \|G\|_{L^{t'}(\mu)} \|r\|_{L^{t}(\mu)}
 \end{aligned}
 \]
where $r\in Rat(K)$ and $C_2$ is a constant. Thus, $B(\lambda _0, \frac{a(\frac{\alpha}{2})\delta}{2}) \subset \Omega .$ This implies $\Gamma _{\frac{\delta}{2}} ^-(w, \frac{\alpha}{4} ) \subset \Omega $ for $\delta$ small enough. Let 
 \[
 \ F(\alpha ) = \{z\in \partial _I^{LG} K\setminus E_1: ~ \exists ~\delta > 0, \text{ such that } \Gamma _{\delta} ^-(z,\alpha) \subset \Omega\},
 \]
then $w\in F(\frac{\alpha}{4})$ and there exists a $\mathcal H^1 |_{\partial _I^{LG} K}$ zero set $E_0$ such that 
 \[
 \ \partial _I^{LG} K\setminus (E_0\cup E_1) \subset F(\frac{\alpha}{4}). 
 \]
It is easy to verify $\mathcal H^1 |_{\partial _I^{LG} K} (F(\alpha _1) \setminus F(\alpha _2) ) = 0$ for $\alpha _1 \ne \alpha _2.$ Let $E_2 = F(\frac{\alpha}{4}) \setminus F(\alpha ),$ then $\mathcal H^1 |_{\partial _I^{LG} K} (E_2) = 0$ and  
 \[
 \ \partial _I^{LG} K\setminus (E_0\cup E_1\cup E_2) \subset \partial _{II}^{LG}.
 \]
The theorem is proved.

The following example is an interesting application of above theorem.

\begin{Example}
A Swiss cheese $K$ can be constructed as
 \[
 \ K = \bar {\mathbb D} \setminus \cup_{n=1}^\infty B(a_n, r_n),
 \] 
where $B(a_n, r_n) \subset \mathbb D,$ $\bar B(a_i, r_i) \cap \bar B(a_j, r_j) = \emptyset $ for $i\ne j,$ $\sum_{n=1}^\infty r_n< \infty ,$ and $K$ has no interior points. Let $\mu$ be the sum of the arc length measures of $\partial \mathbb D$ and  all $\partial B(a_n, r_n).$ Let $\nu$ be the sum of $dz$ on  $\partial \mathbb D$ and all $-dz$ on $\partial B(a_n, r_n).$ For $f\in Rat(K),$ we have
 \[
 \ \int f d\nu = 0.
 \]
Clearly $| \frac{d\nu}{d\mu} | > 0, ~ a.e. ~\mu$ and $\overline{(\frac{d\nu}{d\mu})} \perp R^2(K, \mu),$ so $R^2(K, \mu)$ is irreducible. From Theorem \ref{SBTheorem}, we conclude that
 \[
 \ \underset{\delta\rightarrow 0}{\overline{\lim}}\dfrac{\gamma (\Gamma ^\delta (e^{i\theta}) \setminus K)}{\delta} = 0
 \]
$m$-almost all $e^{i\theta}\in \partial \mathbb D,$ where $\Gamma ^\delta (e^{i\theta})$ is defined in section 2 (right before Theorem \ref{MTheorem1}).
\end{Example}

The example indicates although swiss cheese $K$ has no interior, the portion of $\mathbb D \setminus K$ near $\partial \mathbb D$ is very small.

\begin{Theorem}\label{MTheorem3}
Suppose that $\mu$ is a finite positive measure supported in $K$ and is such that
$abpe(R^t(K,\mu )) = \Omega $ and $R^t(K,\mu )$ is irreducible, where $\Omega$ is a connected region satisfying (3-3). Suppose that upper cone $\Gamma _\delta ^+(w, \alpha )$ is outer for all $w\in \partial _I^{LG} K$ and $\mu = \mu_a + \mu_s$ is the Radon-Nikodym decomposition with respect to $\mathcal H^1 |_{\partial _I^{LG} K},$ where $\mu_a = \frac{1}{2\pi} h\mathcal H^1 |_{\partial _I^{LG} K}$ and $\mu_s \perp \mathcal H^1 |_{\partial _I^{LG} K},$  and $\mu _a(\partial _{II}^{LG} K ) > 0.$ Then:
\newline
(a) If $f \in R^t(K,\mu )$ then the nontangential limit $f^*(z)$ of $f$ exists for $\mu _a |_{\partial _{II}^{LG} K }$-
almost all z, and $f^* = f |_{\partial _{II}^{LG} K }$ as elements of $L^t(\mu |_{\partial _{II}^{LG} K }).$
\newline
(b) Every nonzero rationally invariant subspace $M$ of $R^t(K,\mu )$ has index 1, that is, if $\lambda _0 \in \Omega,$ then $\dim (M/(S_\mu - \lambda _0)M) = 1.$
\newline
If the diameters of all components of $\mathbb C \setminus K$ are bounded away from zero, then by Theorem \ref{SBTheorem}, the above ${\partial _{II}^{LG} K }$ can be replaced by $\partial _{I}^{LG} K .$ 
\end{Theorem}

{\bf Proof:} The proof is the same as in Theorem \ref{MTheorem1} if we apply Proposition  \ref{MProposition1Gen} instead of Proposition \ref{MProposition1}.

The following lemma is an easy exercise.

\begin{Lemma}\label{LemmaEasy}
Let $B(\lambda, \epsilon \delta ) \subset \Gamma _\delta ^-(w, \alpha )$ (or $\Gamma _\delta ^+(w, \alpha )$). Then there are constants $c(\alpha ),~ C(\alpha ) > 0$ that only depend on $\alpha$ and $\|A'\|_\infty$ such that
 \[
 \ \min (\epsilon, c(\alpha ))(\delta + |Re(z - w)|) \le |z - \lambda | \le C(\alpha )(\delta + |Re(z-w)|)
 \]  
for $z\in LG.$
\end{Lemma}

{\bf Proof:} In fact, $C(\alpha ) =1 +  \sqrt{1+ \|A'\|_\infty^2 }$ and $c(\alpha ) =\frac{1 - \alpha \|A'\|_\infty}{\sqrt{1+ \|A'\|_\infty^2 }\sqrt{1+ \alpha^2 }}.$ We leave the details to the reader.

Because we do not have an analogous testing function (such as $(1-\bar\lambda _0 z)^{-\frac{2}{t}}$ in Proposition \ref{MProposition2}) in general, we are not able to get an estimation of the Cauchy transform as in Proposition \ref{MProposition2}. However, our following proposition is enough for us to estimate an upper bound as in (1-2) ((1.4) in \cite{ars}). 
We define a set 
 \[
 \ B\Gamma _\delta ^-(w, \alpha ) = \underset{\lambda _0 \in\Gamma _\delta ^-(w, \frac{\alpha }{2}) \cap (\partial B(w, \frac{\delta}{2}))}{\cup}B(\lambda _0, a(\alpha)\delta)
 \]
where $a(\alpha)$ is defined as in (3-8).

\begin{Proposition}\label{MProposition2Gen} Let $\mu$ be a finite complex  measure with support in $K.$  Suppose that $\mu = \mu_a + \mu_s$ is the Radon-Nikodym decomposition with respect to $\mathcal H^1 |_{\partial _I^{LG} K},$ where $\mu_a = \frac{1}{2\pi} h\mathcal H^1 |_{\partial _I^{LG} K}$ and $\mu_s \perp \mathcal H^1 |_{\partial _I^{LG} K}.$ Suppose upper cone $\Gamma _\delta ^+(w, \alpha )$ is outer for $w\in \partial _I^{LG} K.$ Let $1 < p <\infty, ~ q = \frac{p}{p-1}, ~ f\in C(K), ~ g \in L^q (\mu ),$ and $fg\mu \perp Rat(K).$ Then for $b > 0,$ and $\mathcal H^1 |_{\partial _I^{LG} K}$-almost all $w\in \partial _I^{LG} K,$ there exist $\delta _w >0,$ $E_{\delta}^f \subset B(w, \delta),$ and $\epsilon (\delta ) > 0,$ where $0 < \delta < \delta _w,$ such that $\lim_{\delta\rightarrow 0}\epsilon (\delta ) = 0,$ $\gamma(E_\delta ^f) < \epsilon (\delta ) \delta ,$ and for $\lambda _\delta$ as in (3-4),
\[
\ \begin{aligned}
\ &\left  |\mathcal C \left (fg\mu \right )(\lambda) \right | \le b \delta ^{-\frac{1}{p}} \|f\|_{L^p(\mu )} \\
 + &\dfrac{2 C(\alpha )^{\frac{2}{q}} \delta ^{-\frac{1}{p}} \|f\|_{L^p(\mu )}}{\epsilon_w ^\alpha c_0(\alpha )} \left ( \int \dfrac{\delta}{|Re(z-w) - (\lambda _\delta -w)|^2} |g|^qd\mu _a \right )^{\frac{1}{q}}  
\ \end{aligned}
\]
for all $\lambda\in (B\Gamma _\delta ^-(w, \alpha ) \setminus E_\delta^f )\cap U(g\mu ) ,$ where $\epsilon_w ^\alpha = \min(\epsilon_w, c(\alpha ))$ and $c_0(\alpha ) = \min(a(\alpha ),c(\alpha )),$ where $\epsilon_w$ is as in (3-4), $a(\alpha )$ is from (3-8),  and $c(\alpha ),~ C(\alpha )$ are from Lemma \ref{LemmaEasy}.
\end{Proposition}

{\bf Proof:} Using Lemma \ref{LemmaEasy}, we have the following calculation:
 \[
 \ \begin{aligned}
\ &\left  |\mathcal C (fg\mu _a )(\lambda) - \mathcal C (fg\mu _a )(\lambda _\delta) \right | \\
 \ \le & \int \dfrac{|\lambda - \lambda _\delta|}{|z - \lambda ||z - \lambda _\delta|} |fg|d\mu _a \\
\ \le & \dfrac{2\delta}{\epsilon_w ^\alpha c_0(\alpha )}\int \dfrac{1}{(|Re(z - w)| + \delta )^2} |fg|d\mu _a \\
\ \le & \dfrac{2\delta ^{-\frac{1}{p}}}{\epsilon_w ^\alpha c_0(\alpha )}\int \dfrac{\delta ^\frac{1}{q}}{(|Re(z - w)| + \delta )^\frac{2}{q}} |fg|d\mu _a \\ 
\ \le & \dfrac{2\delta ^{-\frac{1}{p}}\|f\|_{L^p(\mu )}}{\epsilon_w ^\alpha c_0(\alpha )}\left (\int \dfrac{\delta}{(|Re(z - w)| + \delta )^2} |g|^qd\mu _a \right )^\frac{1}{q}\\
 \ \le &\dfrac{2C(\alpha )^{\frac{2}{q}}\delta ^{-\frac{1}{p}} \|f\|_{L^p(\mu )}}{\epsilon_w ^\alpha c_0(\alpha )} \left ( \int \dfrac{\delta}{|Re(z-w) - (\lambda _\delta -w)|^2} |g|^qd\mu _a \right )^{\frac{1}{q}},  
\ \end{aligned}
 \]
where the last step also follows Lemma \ref{LemmaEasy}. The rest of proof is the same as in the proof of Proposition \ref{MProposition2}.

\begin{Theorem}\label{MTheorem4} 
Suppose that $\mu$ is a finite positive measure supported in $K$ and is such that
$abpe(R^t(K,\mu )) = \Omega $ and $R^t(K,\mu )$ is irreducible, where $\Omega$ is a connected region satisfying (3-3). Suppose that upper cone $\Gamma _\delta ^+(w, \alpha )$ is outer for all $w\in \partial _I^{LG} K$ and $\mu = \mu_a + \mu_s$ is the Radon-Nikodym decomposition with respect to $\mathcal H^1 |_{\partial _{I}^{LG} K},$ where $\mu_a = \frac{1}{2\pi} h\mathcal H^1 |_{\partial _{I}^{LG} K }$ and $\mu_s \perp \mathcal H^1 |_{\partial _{I}^{LG} K},$  and $\mu _a(\partial _{II}^{LG} K ) > 0.$ Then:
\newline 
(a) For $t = 1,$ there are constants $C(w) > 0$ (depending on $G$) such that  
 \[
 \ \underset{\Gamma ^-(w, \frac{\alpha}{2}) \ni \lambda \rightarrow w}{\overline\lim} |\lambda -w |M_\lambda \le  \dfrac{C(w)}{h(w)}\tag{3-11} 
 \]
for $\mu_a$-almost all $w\in \partial _{II}^{LG} K.$
\newline 
(b) For $t > 1,$ there are constants $C_0(\alpha ) > 0$ (depending on $\alpha $ and $\|A'\|_\infty$) such that  
 \[
 \ \underset{\Gamma ^-(w, \frac{\alpha}{2})\ni \lambda \rightarrow w}{\overline\lim} |\lambda -w |^{\frac{1}{t}} M_\lambda \le  \dfrac{C_0(\alpha )/(\epsilon _w \epsilon _w^\alpha )}{h(w)^{\frac{1}{t}}}\tag{3-12} 
 \]
for $\mu_a$-almost all $w\in \partial _{II}^{LG} K,$ where $\epsilon _w$ is as in  (3-4) and $\epsilon _w^\alpha$ is from Proposition \ref{MProposition2Gen}.
\newline 
If the diameters of all components of $\mathbb C \setminus K$ are bounded away from zero, then by Theorem \ref{SBTheorem}, the above $\partial _{II}^{LG} K$ can be replaced by $\partial _{I}^{LG} K.$
\end{Theorem}

{\bf Proof:} (a) Let $w\in \partial _{II}^{LG} K$ so that we can apply Proposition  \ref{MProposition1Gen} for $w$ and $b = \frac{ |G(w)h(w)|}{2}.$ There exist $\delta _w > 0,$ $E_{\delta}\subset B(w, \delta),$ and $\epsilon (\delta ) > 0,$ where $0 < \delta < \delta _w,$  such that $\lim_{\delta\rightarrow 0}\epsilon (\delta ) = 0,$ $\gamma(E_\delta) < \epsilon (\delta ) \delta ,$ and
\[
\ \left  |\mathcal C(G\mu ) (\lambda) \right | \ge \frac{|G(w)h(w)|}{2} \tag{3-13}
\]
for all $\lambda\in (\Gamma _\delta ^-(w, \alpha )\setminus E_\delta )\cap U(G\mu ) ,$ where $\Gamma _\delta ^-(w, \alpha ) \subset \Omega.$
Now choose $\delta$ to be small enough so that $\epsilon (\delta ) < a(\alpha )\epsilon _1,$ where $\epsilon _1$ is as in Lemma \ref{lemmaARS} and $a(\alpha )$ is from (3-8). Let $\lambda _0 \in \Gamma _{\delta} ^-(w, \frac{\alpha}{2} )$ and $|\lambda _0 - w| = \frac{\delta}{2},$ then $B(\lambda _0, a(\alpha )\delta) \subset \Gamma _{\delta} ^-(w, \alpha) \subset  \Omega,$ where $\delta$ is small enough. Since $\gamma (B(\lambda _0, a(\alpha )\delta) \cap E_\delta ) < \epsilon_1 a(\alpha )\delta, $ from Lemma \ref{lemmaARS}, (2-16),  and (3-13). we conclude that for $\lambda \in B(\lambda _0, \frac{a(\alpha )\delta}{2})$ and $r\in Rat(K),$  we have
\[ 
 \ \begin{aligned}
 \  |r (\lambda) | \le & \dfrac{C_1}{(a(\alpha )\delta) ^2} \int _{B (\lambda _0, a(\alpha )\delta ) \setminus E_\delta} |r(z)| \dfrac{dA(z)}{\pi} \\
 \ \le & \dfrac{C_1}{\pi a(\alpha )^2\delta ^2} \int _{B (\lambda _0, a(\alpha )\delta ) \setminus E_\delta} \dfrac{|\mathcal C(rG\mu ) (z)|}{|\mathcal C(G\mu ) (z)|} dA(z)\\
 \ \le & \dfrac{2C_1}{\pi |G(w)h(w)|a(\alpha )^2\delta ^2} \int \int _{B (\lambda _0, a(\alpha )\delta )} \dfrac{1}{|z-u|} dA(z) |r(u)||G(u)|d\mu (u)\\
 \ \le & \dfrac{C_2}{|G(w)h(w)|\delta } \int |r(u)||G(u)|d\mu (u),
 \ \end{aligned}
 \]
where $C_1,$ $C_2,$ $C_3,...$ stand for absolute constants, and hence,
 \[
 \ |\lambda - w|  |r (\lambda) | \le \dfrac{C_3}{|G(w)h(w)|} \int |r(u)||G(u)|d\mu (u) 
	\]
for $\lambda  \in \Gamma _{\delta} ^-(w, \frac{\alpha}{2} )$ and $|\lambda  - w| = \frac{\delta}{2}.$	
	Let $C(w) = \dfrac{C_3\|G\|_{L^\infty (\mu )}}{|G(w)|},$ we get
	 \[
 \ \underset{\Gamma ^-(w, \frac{\alpha}{2})\ni \lambda \rightarrow w}{\overline{\lim}}|\lambda - w|M_\lambda  \le  \dfrac{C(w)}{h(w)} .
	\]
	
	(b) By Proposition \ref{MProposition1Gen} and Proposition \ref{MProposition2Gen}, for $b > 0,$ and $\mathcal H^1 |_{\partial _I^{LG} K}$-almost all $w\in \partial _I^{LG} K,$ there exist $\delta _w >0,$ $E_{\delta} \subset B(w, \delta),$ $E_{\delta}^r \subset B(w, \delta),$ and $\epsilon (\delta ) > 0,$ where $0 < \delta < \delta _w,$ such that $\lim_{\delta\rightarrow 0}\epsilon (\delta ) = 0,$ $\gamma(E_\delta) < \epsilon (\delta ) \delta ,$ $\gamma(E_\delta ^r) < \epsilon (\delta ) \delta ,$ 
	\[
\  |\mathcal C(G\mu ) (\lambda) - L(w)G(w)h(w)| \le b,\tag{3-14}
\]
for all $\lambda\in (\Gamma _\delta ^-(w, \alpha )\setminus E_\delta )\cap U(G\mu ) ,$
	and for $\lambda _\delta$ as in (3-4),
\[
\ \begin{aligned}
\ &\left  |\mathcal C \left (rg\mu \right )(\lambda) \right | \le b \delta ^{-\frac{1}{t}} \|r\|_{L^t(\mu )} \\
 + &\dfrac{2C(\alpha )^{\frac{2}{t'}}\delta ^{-\frac{1}{t}} \|r\|_{L^t(\mu )}}{\epsilon_w ^\alpha c_0(\alpha )} \left ( \int \dfrac{\delta}{|Re(z-w) - (\lambda _\delta -w)|^2} |G|^{t'}d\mu _a \right )^{\frac{1}{t'}}  
\ \end{aligned} \tag{3-15}
\]
for all $\lambda\in (B\Gamma _\delta ^-(w, \alpha ) \setminus E_\delta^r )\cap U(|G|^{t'}\mu ) ,$
	
From Plemelj's formula (3-2), we have the following calculation:
\[
\ \begin{aligned}
\ & \underset{\delta\rightarrow 0}{\overline\lim}\int \dfrac{\delta}{|Re(z-w) - (\lambda _\delta -w)|^2} |G|^{t'}d\mu _a \\
\ = & \underset{\delta\rightarrow 0}{\overline\lim}\dfrac{i\delta}{2 Im(\lambda _\delta)} (\dfrac{1}{2\pi i} \mathcal C (|G|^{t'}h\sqrt{1 + (A'(x))^2}dx) (\lambda _\delta - w - Re(w)) \\
\ &- \dfrac{1}{2\pi i} \mathcal C (|G|^{t'}h\sqrt{1 + (A'(x))^2}dx) (\bar\lambda _\delta - \bar w  - Re(w) )) \\
\ \le &\dfrac{|G(w)|^{t'} h(w)\sqrt{1 + (A'(Re(w)))^2}}{2\epsilon_w}.    
\ \end{aligned}\tag{3-16}
\]
Therefore, for $\eta > 0,$ if $\delta$ is small enough, we conclude
 \[
 \  \int \dfrac{\delta}{|Re(z-w) - (\lambda _\delta -w)|^2} |G|^{t'}d\mu _a < \dfrac{|G(w)|^{t'} h(w)\sqrt{1 + (A'(Re(w)))^2}}{2\epsilon_w} + \eta .\tag{3-17}
 \]
Combining (3-14), (3-15), and (3-17), for $\delta$ small enough, $\lambda_0\in (\partial B(w,\frac{\delta}{2}))\cap \Gamma _\delta ^-(w, \frac{\alpha }{2}) ,$ $B(\lambda_0, a(\alpha ) \delta ) \subset \Gamma _\delta ^-(w, \alpha ), $ and $\lambda \in (B(\lambda_0, a(\alpha ) \delta ) \setminus (E_\delta \cup E_\delta^r))\cap U(|G|^{t'}\mu),$ we get
 \[
 \ \dfrac{|r(\lambda )|}{\|r\|_{L^t(\mu )}} \le \dfrac{b \delta ^{-\frac{1}{t}} + \dfrac{2C(\alpha )^{\frac{2}{t'}}\delta ^{-\frac{1}{t}}}{\epsilon_w ^\alpha c_0(\alpha )} \left ( \dfrac{|G(w)|^{t'} h(w)\sqrt{1 + \|A'\|_\infty^2}}{2\epsilon_w} + \eta \right)^{\frac{1}{t'}}}{|G(w)| h(w) -b }.\tag{3-18}
 \]
From semi-additivity of (2-2), we see
 \[
 \ \gamma (E_\delta \cup E_\delta^r) \le A_T(\gamma (E_\delta ) + \gamma ( E_\delta^r)) \le 2A_T \epsilon (\delta ) \delta .
 \] 
Let $\delta $ be small enough so that $\epsilon (\delta ) < \frac{a(\alpha)}{2A_T}\epsilon _1,$ where $\epsilon _1$ is as in Corollary \ref{CorollaryDSet}. From Corollary \ref{CorollaryDSet}, we conclude that (3-18) holds for all $\lambda \in (B(\lambda_0, \frac{a(\alpha ) \delta}{2}).$ Hence, for $\delta$ small enough, $\lambda\in (\partial B(w,\frac{\delta}{2}))\cap \Gamma _\delta ^-(w, \frac{\alpha }{2}),$
 \[
 \ |\lambda - w|^{\frac{1}{t}} M_ \lambda \le \dfrac{b  + \dfrac{2C(\alpha )^{\frac{2}{t'}}}{\epsilon_w ^\alpha c_0(\alpha )} \left ( \dfrac{|G(w)|^{t'} h(w)\sqrt{1 + (A'(Re(w)))^2}}{2\epsilon_w} + \eta \right)^{\frac{1}{t'}}}{|G(w)| h(w) -b }
 \]
Therefore, there exists a constant $C_0(\alpha ) > 0$ that only depends on $\alpha$ and $\|A'\|_\infty$ so that 
	\[
 \ \underset{ \Gamma ^-(w,\frac{\delta}{2})\ni \lambda\rightarrow w}{\overline{\lim}} |\lambda - w|^{\frac{1}{t}} M_ \lambda \le \dfrac{C_0(\alpha )}{\epsilon_w \epsilon_w  ^\alpha h(w)^\frac{1}{t}}
	\]
	for $\mathcal H^1 |_{LG}$-almost all $w \in \partial_{II}^{LG} K.$

	For the lower bound, we do have testing functions $f_1^\delta (z) = (z - \lambda _\delta ) ^{-2}\in R^1(K,\mu)$ and $f_2^\delta (z) = (z - \lambda _\delta ) ^{-1}\in R^2(K,\mu).$ The following proposition estimates their norms.

\begin{Proposition}\label{LemmaKTGen}
Let $\mu$ be a finite positive measure with support in $K.$
Suppose that $\mu = \mu_a + \mu_s$ is the Radon-Nikodym decomposition with respect to $\mathcal H^1 |_{\partial _I^{LG} K},$ where $\mu_a = \frac{1}{2\pi} h\mathcal H^1 |_{\partial _I^{LG} K}$ and $\mu_s \perp \mathcal H^1 |_{\partial _I^{LG} K},$  and $\mu _a(\partial _{II}^{LG} K ) > 0.$ Suppose that $\Gamma _\delta ^+(w,\alpha )$ is outer for $w\in\partial _{II}^{LG} K,$
then there exists a constant $C_1(\alpha ) > 0$ that only depends on $\alpha$ and $\|A'\|_\infty$ such that
 \[
 \ \underset{\delta \rightarrow 0}{\overline\lim} \int \dfrac{\delta}{|z - \lambda _\delta|^2} d \mu \le \dfrac{C_1(\alpha )}{\epsilon_w(\epsilon_w^\alpha)^2}h(w).
 \]
for $\mu _a$-almost all $w\in \partial _{II}^{LG} K,$ where $\lambda _\delta$ and $\epsilon_w$ are from (3-4).
\end{Proposition}

{\bf Proof:} The proposition follows from Corollary \ref{CorollaryKTGen}, Lemma \ref{LemmaEasy}, and the same proof of (3-16). 

So we have lower bounds for $R^1(K,\mu)$ and $R^2(K,\mu)$ as the following:

For $t = 1,$
 \[
 \ \underset{\Gamma ^-(w, \frac{\alpha}{2}) \ni \lambda \rightarrow w}{\underline\lim} |\lambda -w | M_\lambda \ge \underset{\delta \rightarrow 0}{\underline\lim} \dfrac{|f_1^\delta (\lambda )|}{\|f_1^\delta\|_{L^1(\mu)}} \ge \dfrac{\epsilon_w(\epsilon_w^\alpha)^2}{4C_1(\alpha )h(w)}.
 \]

For $t = 2,$
\[
 \ \underset{\Gamma ^-(w, \frac{\alpha}{2}) \ni \lambda \rightarrow w}{\underline\lim} |\lambda -w |^\frac{1}{2} M_\lambda \ge \underset{\delta \rightarrow 0}{\underline\lim} \dfrac{|f_2^\delta (\lambda )|}{\|f_2^\delta\|_{L^2(\mu)}} \ge \dfrac{\sqrt{\epsilon_w}\epsilon_w^\alpha}{2\sqrt{C_1(\alpha )h(w)}}.
 \]
For $t \ne 1$ and $t \ne 2,$ if $w$ is a boundary point of $\mathbb C\setminus K,$ then we can define a similar testing function and have corresponding lower bounds. However, if $w$ is an inner boundary point, we do not have such a testing function to estimate the lower bounds.

\bibliography{Bibliography}

\end{document}